\documentclass[12pt,righttagstyle]{amsart}

\usepackage{amssymb,latexsym}

\newcommand{\comment}[1]{}

\newcommand{\epf}{ $\Box$\medskip}

\newtheorem{thm}{Theorem}[section]
\newtheorem{prop}[thm]{Proposition}

\newtheorem{lem}[thm]{Lemma}
\newtheorem{defn}[thm]{Definition}

\newtheorem{notation}[thm]{Notation}
\begin{document}

\title[Integrable fractional mean functions]{Integrable fractional mean functions on spaces of homogeneous type.}

\author[J. Feuto]{Justin Feuto}
\address{UFR DE MATHEMATIQUES ET INFORMATIQUE,
 Universit\'e de Cocody, 22BP 1194 Abidjan 22, R\'epublique de C\^ote d'Ivoire}
\email{{\tt justfeuto@yahoo.fr}}
\author[I. Fofana]{Ibrahim Fofana}
\address{UFR DE MATHEMATIQUES ET INFORMATIQUE,
 Universit\'e de Cocody, 22BP582 Abidjan 22, R\'epublique de C\^ote d'Ivoire}
\email{{\tt fofana\_ib\_math\_ab@yahoo.fr}}
\author[K. Konin]{Konin Koua}
\address{UFR DE MATHEMATIQUES ET INFORMATIQUE,
 Universit\'e de Cocody, 22BP582 Abidjan 22, R\'epublique de C\^ote d'Ivoire}
\email{{\tt kroubla@yahoo.fr}}

\thanks{}

\begin{abstract}
The class of Banach spaces $\left( L^{q},L^{p}\right) ^{\alpha }(X,d,\mu )$, 
$1\leq q\leq \alpha \leq p\leq \infty ,$ introduced in \cite{F1} in
connection with the study of the continuity of the fractional maximal
operator of Hardy-Littlewood and of the Fourier transformation in the case $%
X=\mathbb{R}^{n}$ and $\mu $ is the Lebesgue measure, was generalized in 
\cite{FFK} to the setting of homogeneous groups. We generalize it here to spaces of homogeneous type and we prove that the results obtained in \cite{FFK} such as relations between these spaces and Lebesgue spaces, weak Lebesgue and Morrey spaces, remain true.  

\medskip

{\scshape R\'esum\'e.} La classe $\left( L^{q},L^{p}\right) ^{\alpha }(X,d,\mu )$, 
$1\leq q\leq \alpha \leq p\leq \infty ,$ introduite dans \cite{F1} en liaison avec l'\'etude de la continuit\'e de l'op\'erateur maximal fractionnaire de Hardy-Littlewood et de la transformation de Fourier dans le cas o\`u $X=\mathbb R^{n}$ et $\mu$ la mesure de Lebesgue, a \'et\'e g\'en\'eralis\'ee dans \cite{FFK} au cas o\`u $X$ est un groupe homog\`ene. Nous le g\'en\'eralisons ici aux espaces de type homog\`ene et nous montrons que la plupart des r\'esultats obtenus dans \cite{FFK}, tels que les liens entre ces espaces et les espaces de Lebesgue, de Lebesgue faible et de Morrey, demeurent valides.
\end{abstract}

\maketitle

\section{Introduction}
In \cite{M}, Muckenhoupt raised the problem of
characterizing weight functions $u$ and $v$ for which the inequality

\begin{equation}
\int_{-\infty}^{+\infty }\left\vert \widehat{f}(x)\right\vert ^{p}u(x)dx\leq
C\int_{-\infty}^{+\infty }\left\vert f(x)\right\vert ^{p}v(x)dx\label{mf}
\end{equation}
holds for every $f$ in the Lebesgue space $L^{p}(\mathbb{R}).$ 

Aguilera and Harboure showed in \cite{AH} that, in the case $v=1$ and $1<p<2$, a necessary condition for (\ref{mf}) is 

\begin{equation}
\left[ \underset{k=-\infty }{\overset{k=+\infty }{\sum }}\left(
\int_{rk}^{r(k+1)}u(x)\right) ^{b}\right] ^{\frac{1}{b}}\leq Cr^{p-1},\ \ \ r>0
\label{0.0001}
\end{equation}

where $b=\frac{2}{2-p}.$

Let us assume that $n$ is a positive integer and $1\leq q\leq\alpha\leq p\leq\infty$. For any Lebesgue-measurable function $f$ on $\mathbb R^{n}$, we set
\begin{equation}
\left\| f\right\| _{q,p,\alpha }=\left\{ 
\begin{array}{lll}
\underset{r>0}{\sup }r^{n\left( \frac{1}{\alpha }-\frac{1}{q}\right) }\left[ 
\underset{k\in \mathbb{Z}^{n}}{\sum }\left( \left\| f\chi
_{_{I_{k}^{r}}}\right\| _{q}\right) ^{p}\right] ^{\frac{1}{p}} & \text{if }
& p< \infty \\ 
\underset{r>0}{\sup }r^{n\left( \frac{1}{\alpha }-\frac{1}{q}\right) }%
\underset{k\in \mathbb{Z}^{n}}{\sup }\left\| f\chi _{_{I_{k}^{r}}}\right\|
_{q} & \text{if} & p=\infty%
\end{array}
\right..\label{normeuclidien}
\end{equation}
Here $I_{k}^{r}=\overset{n}{\underset{j=1}{\Pi }}\left[ k_{j}r\ ,\ \left(
k_{j}+1\right) r\right) $, $k=\left(k_{1},\ldots,k_{n}\right) \in \mathbb{Z}^{n}$,  $r>0$ and $\left\|\cdot\right\|_{q}$ denotes the usual norm on the Lebesgue space $L^{q}(\mathbb R^{n})$. We denote by  $L_{0}(\mathbb{R}^{n})$ the complex vector space of equivalent classes (modulo equality Lebesgue almost everywhere) of Lebesgue measurable complex-valued functions on $\mathbb{R}^{n}.$ It is clear that $\left\|\cdot\right\|_{q,p,\alpha}$ may be looked at as a map of $L_{0}(\mathbb R^{n})$ into $\left[0,\infty\right]$. We define 
\begin{equation}
\left( L^{q},L^{p}\right) ^{\alpha }\left( \mathbb{R}^{n}\right) =\left\{
f\in L_{0}(\mathbb{R}^{n})/\left\| f\right\| _{q,p,\alpha }<\infty \right\}.
\end{equation}

Fofana has proved in \cite{F1} that ($\left( L^{q},L^{p}\right) ^{\alpha }\left(\mathbb{R}^{n}\right) ,\left\|\cdot\right\| _{q,p,\alpha })$ is a complex Banach space and that the Lebesgue spaces $L^{\alpha }\left( \mathbb{R}^{n}\right) $, the Morrey spaces $M_{q}^{n\left( 1-\frac{1}{\alpha }\right) }\left(\mathbb{R}^{n}\right)$  and the Lorenz spaces $L^{\alpha ,\infty }(\mathbb{R}^{n})$ (in the case $q<\alpha<p$) are its sub-spaces.

Note that condition (\ref{0.0001}) can be written as $u\in\left( L^{1},L^{b}\right) ^{\frac{1}{2-p}},$ with $b=\frac{2}{2-p}$.

Further results on Fourier transform may be expressed in the setting of $(L^{q},L^{p})^{\alpha}(\mathbb R^{n})$ and related spaces of Radon measures (see \cite{F3}, \cite{L}). These spaces are also related to $L^{q}-L^{p}$ multiplier problems (see \cite{FKK}, \cite{O}) and well-suited to establish norm inequalities for fractional maximal functions \cite{F4}.

It is clear that $(L^{q},\ell^{p})^{\alpha}(\mathbb R^{n})$ is a subspace of the so-called amalgam space of Wiener $(L^{q},\ell^{p})(\mathbb R^{n})$, defined by 
\begin{equation}
(L^{q},\ell^{p})^{\alpha}(\mathbb R^{n})=\left\{f\in L_{0}(\mathbb R^{n}):\  _{1}\left\|f\right\|_{q,p}<\infty\right\}
\end{equation}
where for $r>0$
\begin{equation}
_{r}\left\| f\right\| _{q,p}=\left\{ 
\begin{array}{lll}
\left[\sum_{n\in\mathbb Z^{n}}\left(\left\|f\chi_{I^{r}_{k}}\right\|_{q}\right)^{p}\right]^{\frac{1}{p}} & \text{ if }  &  p<\infty \\ 
\underset{k\in\mathbb Z^{n}}{\sup}\left\|f\chi_{I^{r}_{k}}\right\|_{q} & \text{ if }& p=\infty%
\end{array}
\right..\label{namalgamR}
\end{equation}

These amalgam spaces have been used by Wiener (see \cite{Wi}) in connection with Tauberian theorems. Long after, Holland undertook their systematic study (see \cite{H}). Since then, they have been extensively studied (see the survey paper \cite{FS} and the references therein) and generalized to locally compact groups (see \cite{FFK}, \cite{BS}, \cite{BD}). They may be looked at as spaces of functions that behave locally as elements of $L^{q}(\mathbb R^{n})$ and globally as belonging to $L^{p}(\mathbb R^{n})$. Taking this into  account, Feichtinger has introduced Banach spaces whose elements belong locally to some Banach space, and globally to another (see \cite{Fe}).

Replacing $\mathbb R^{n}$ by a homogeneous group $G$, the authors have defined and studied $\left( L^{q},L^{p}\right) ^{\alpha }(G)$ spaces (see \cite{FFK}). They proved that results obtained in \cite{F1} remain valid for $\left( L^{q},L^{p}\right) ^{\alpha }(G)$.

In the present paper, we extend the definition of these spaces to spaces of homogeneous type. In this setting, we obtain interesting links between $\left( L^{q},L^{p}\right) ^{\alpha }(X)$ and classical Banach function spaces. 

These spaces are well suited for studying norm inequalities on fractional maximal operators. Actually, in \cite{F-F-K} we established some continuity properties for these operators from $(L^{q},L^{p})^{\alpha}(X)$ to weak-Lebesgue spaces, which extended in this context analogous results known in the Euclidean case (see \cite{F4}, \cite{MW}).

The paper is organized as follows. Section 2 contains definitions and the main results, whose proofs are given in Section 4. Section 3 is devoted to  auxiliary results. 

 Throughout the paper, $C$ denotes positive constants that are independent of the main parameters involved, with values which may differ from line to
line. Constants with subscripts, such as $C_{1}$, do not change in different occurrences. 

\section{Definitions-Results}
A space of homogeneous type $(X,d,\mu)$ is a quasi metric space $(X,d)$ endowed with a non negative Borel measure $\mu$ satisfying the doubling condition 
\begin{equation}
0<\mu \left( B_{\left(x,2r\right) }\right) \leq C\mu \left(B_{\left(
x,r\right)}\right) <\infty ,\ x\in X\text{ and } r>0, \label{0.2}
\end{equation}
 where $B_{\left( x,r\right) }=\left\{y\in X:d\left( x,y\right) <r\right\} $ is the ball with center $x$ and radius $r>0$. 
If $C'_{\mu}$ is the smallest constant $C$ for which (\ref{0.2}) 
holds, then $D_{\mu}=\log _{2}C'_{\mu}$ is called the
doubling order of $\mu$. It is known (see \cite{sw}) that  for all balls $B_{2}\subset B_{1}$ 
\begin{equation}
\frac{\mu \left( B_{1}\right) }{\mu \left( B_{2}\right) }\leq C_{\mu}\left( \frac{r\left(
B_{1}\right) }{r\left( B_{2}\right) }\right) ^{D_{\mu}},
\label{0.05}
\end{equation}
where $r(B)$ denote the radius of the ball $B$ and $C_{\mu}=C'_{\mu}(2\kappa)^{D_{\mu}}$, $\kappa \geq 1$ being a constant such that 
\begin{equation}
d(x,y)\leq \kappa \left( d(x,z)+d(z,y)\right),\ \ \ x,y,z\in X.  \label{0.001}
\end{equation} 

 Two quasi metrics $d$ and $\delta$ on $X$ are said to be equivalent if there exists constants $C_{1}>0$ and $C_{2}>0$ such that
\begin{equation*}
C_{1}d(x,y)\leq\delta(x,y)\leq C_{2}d(x,y),\ \ \ x,y\in X.
\end{equation*}

Observe that topologies defined by equivalent quasi metrics on $X$ are equivalent. It is shown in \cite{MS} that on any space of homogeneous type $(X,d,\mu)$, there is a quasi metric $\delta$ equivalent to $d$ for which balls are open sets.
\medskip

In the sequel we assume that $X=\left( X,d,\mu \right) $ is a fixed space of homogeneous type and:\\
$\bullet$ all balls $B_{\left( x,r\right) }=\left\{y\in X:d\left( x,y\right) <r\right\} $ are open subsets of $X$ endowed with the $d$-topology,\\
$\bullet$ $(X,d)$ is separable,\\ 
$\bullet \ \mu(X)=\infty,$\\
$\bullet \  B_{(x,R)}\setminus B_{(x,r)}\neq \emptyset,\ 0<r<R<\infty$ and $x\in X$.

As proved in \cite{W}, the last assumption implies that there exists two constants $\tilde{C}_{\mu}>0$ and $\delta_{\mu}>0$ such that for all balls $B_{2}\subset B_{1}$ of $X$
 \begin{equation}
 \frac{\mu(B_{1})}{\mu(B_{2})}\geq \tilde{C}_{\mu}\left(\frac{r(B_{1})}{r(B_{2})}\right)^{\delta_{\mu}}.\label{revd}
 \end{equation}

For $1\leq p\leq \infty $, $\left\|\cdot\right\| _{p}$\ denotes the usual norm on the Lebesgue space $L^{p}(X)$. 

For any $\mu$-measurable function $f$ on $X$, we set:\\
$\bullet\ \lambda_{f}(\alpha)=\mu\left(\left\{x\in X:\left|f(x)\right|>\alpha\right\}\right),\ \alpha>0$,\\
$\bullet\ f_{\ast}(t)=\inf\left\{\alpha>0:\lambda_{f}(\alpha)\leq t\right\},\ t>0$,\\
$\bullet\ f^{\ast}(t)=\frac{1}{t}\int^{t}_{0}f_{\ast}(u)du,\ t>0$,\\
$\bullet\left\|f\right\|_{p,q}=\left\{
\begin{array}{ccc}
\left[\frac{p}{q}\int^{\infty}_{0}\left(t^{\frac{1}{p}}f^{\ast}(t)\right)^{q}\frac{dt}{t}\right]^{\frac{1}{q}} & \text{ if } & 1\leq p,q<\infty\\
\sup_{t>0}t^{\frac{1}{p}}f^{\ast}(t) & \text{ if } & 1\leq p\leq\infty \text{ and } q=\infty
\end{array}
\right.$.

 Let $L_{0}(X)$ be the complex vector space of equivalent classes (modulo equality $\mu -$almost everywhere) of $\mu -$measurable complex-valued functions on $X$. Then $\left\|\cdot\right\|_{p,q}$ is a map from $L_{0}(X)$ into $\left[0,\infty\right]$. It is known (see \cite{SWe}) that: \\
$\bullet$ for $1<p,q\leq\infty$, the space $L^{p,q}(X)=\left\{f\in L_{0}(X):\left\|f\right\|_{p,q}<\infty\right\}$ endowed with $f\mapsto\left\|f\right\|_{p,q}$, is a complex Banach space (called Lorentz space),\\
$\bullet \ f\mapsto\left\|f\right\|^{\ast}_{p,q}=\left\{\begin{array}{ccc}
\left[\frac{p}{q}\int^{\infty}_{0}\left(t^{\frac{1}{p}}f_{\ast}(t)\right)^{q}\frac{dt}{t}\right]^{\frac{1}{q}} & \text{ if } & 1\leq p,q<\infty\\
\sup_{t>0}t^{\frac{1}{p}}f_{\ast}(t) & \text{ if } & 1\leq p\leq\infty \text{ and } q=\infty
\end{array}
\right.$\\
is a quasi-norm on $L^{p,q}(X)$ equivalent to $\left\|\cdot\right\|_{p,q}$,\\
$\bullet\ \underset{t>0}{\sup}t^{\frac{1}{p}}f_{\ast}(t)=\underset{\alpha>0}{\sup}\alpha\lambda_{f}(\alpha)^{\frac{1}{p}}$.
 
In the sequel we assume that $1\leq q\leq\alpha\leq p\leq\infty$.
\begin{notation}
For any $\mu$-measurable function $f$ on $X$ and any $r>0$, we put
\begin{equation}
 _{r}\left\| f\right\| _{q,p,\alpha}=\left\{ 
\begin{array}{lll}
\left[ \int_{X}\left(\mu(B_{(y,r)})^{\frac{1}{\alpha}-\frac{1}{p}-\frac{1}{q}} \left\| f\chi _{_{B_{\left( y,r\right)
}}}\right\| _{q}\right) ^{p}d\mu (y)\right] ^{\frac{1}{p}} & \text{if} & 
p<\infty \\ 
\underset{y\in X}{\sup}\mu(B_{(y,r)})^{\frac{1}{\alpha}-\frac{1}{q}}\left\| f\chi _{_{B_{\left( y,r\right)
}}}\right\| _{q} & \text{if} & p=\infty%
\end{array}
\right. ,\label{namalgamX}
\end{equation}
where $\chi _{_{B_{\left( y,r\right)}}}$ denotes the characteristic function of $B_{\left( y,r\right)}$, and we use the convention $\frac{1}{q}=0$ if $q=\infty$.
\end{notation}
\begin{thm}\label{theorem2.2}
For any $\mu$-measurable function $f$ on $X$ and $r>0$, we have $_{r}\left\| f\right\| _{q,p,\alpha}=0$ if and only if $f=0$ $\mu$-almost everywhere.
\end{thm}
By the previous result we may (and shall) look at $ _{r}\left\| \cdot\right\| _{q,p,\alpha}$ as a map from $L_{0}(X)$ into $\left[0,\infty\right]$.

\begin{notation}
For $r>0$, we define
\begin{equation*}
\left( L^{q},L^{p}\right)^{\alpha} _{r}(X)=\left\{
f\in L_{0}(X):\;_{r}\left\| f\right\| _{q,p,\alpha}<\infty
\right\} .
\end{equation*}
\end{notation}

\begin{thm}\label{eamalgamX}
For any positive real number $r$, $\left(\left(L^{q},L^{p}\right)^{\alpha}_{r}(X),\ _{r}\left\|\cdot\right\|_{q,p,\alpha}\right)$ is a complex Banach space.
\end{thm}

As in the Euclidean case we have the following results.
\begin{thm}\label{comparaisonnormamalgame}
Let $r$ be a positive real number, $1\leq q_{1}<q_{2}\leq\alpha$ and $\alpha\leq p_{1}<p_{2}<\infty$. Then
\begin{equation}
_{r}\left\|\cdot\right\|_{q_{1},p,\alpha}\leq\ _{r}\left\|\cdot\right\|_{q_{2},p,\alpha} \label{9}
\end{equation}
and
\begin{equation}
_{r}\left\|\cdot\right\|_{q,\infty,\alpha}\leq C\  _{r}\left\|\cdot\right\|_{q,p_{2},\alpha}\leq C'\  _{r}\left\|\cdot\right\|_{q,p_{1},\alpha},\label{10}
\end{equation}
where $C>0$ and $C'>0$ are constants independent of $r$.
\end{thm}

\begin{thm}\label{comparaisonnormamalgamelebesgue}
There is a constant $C>0$ such that 
\begin{equation}
_{r}\left\|\cdot\right\|_{q,p,\alpha}\leq C \  _{r}\left\|\cdot\right\|_{\alpha}\label{11}
\end{equation}
for any real number $r>0$.
\end{thm}
$(L^{q},L^{p})^{\alpha}_{r}(X)$ is actually a generalization of the Wiener amalgam space $(L^{q},\ell^{p})(\mathbb R^{n})$. This appears clearly when we compare $_{r}\left\|\cdot\right\|_{q,p}$ as define in (\ref{namalgamR}), to the norm $\left\|\cdot\right\|^{d_{m_{r}}}_{q,p,\alpha}$ which is equivalent to $_{r}\left\|\cdot\right\|_{q,p,\alpha}$ (see Proposition \ref{normequivalente}). Now we define a subspace of $(L^{q},L^{p})^{\alpha}_{r}(X)$ which generalizes $(L^{q},\ell^{p})^{\alpha}(\mathbb R^{n})$.

\begin{defn}\label{espacefractionnaire}
We set 
\begin{equation*}
(L^{q},L^{p})^{\alpha}(X)=\left\{f\in L_{0}(X):\left\|f\right\|_{q,p,\alpha}<\infty\right\},
\end{equation*}
where $\left\|f\right\|_{q,p,\alpha}=\sup_{r>0}\ _{r}\left\|f\right\|_{q,p,\alpha}$.
\end{defn}
From Definition \ref{espacefractionnaire}, Theorem \ref{eamalgamX}, Theorem \ref{comparaisonnormamalgamelebesgue} and Theorem \ref{comparaisonnormamalgame} the following results are straightforward.

\begin{thm}\label{theorem2.8}
\begin{enumerate}
\item[$a)$] $\left((L^{q},L^{p})^{\alpha}(X),\left\|\cdot\right\|_{q,p,\alpha}\right)$ is a complex Banach space and there exists a constant $C>0$ such that
\begin{equation}
\left\|\cdot\right\|_{q,p,\alpha}\leq C\left\|\cdot\right\|_{\alpha}.\label{inclusionlebesgueamalgam}
\end{equation}
\item[$b)$] Assume that $1\leq q_{1}<q_{2}\leq\alpha\leq p_{1}<p_{2}\leq\infty$. Then
\begin{equation*}
\left\|\cdot\right\|_{q_{1},p,\alpha}\leq C\left\|\cdot\right\|_{q_{2},p,\alpha}
\end{equation*}
and
\begin{equation*}
\left\|\cdot\right\|_{q,p_{2},\alpha}\leq C\left\|\cdot\right\|_{q,p_{1},\alpha},
\end{equation*}
for some constant $C>0$.
\end{enumerate}
\end{thm}

The continuous embedding of $L^{\alpha}(X)$ into $(L^{q},L^{p})^{\alpha}(X)$ expressed by Inequality (\ref{inclusionlebesgueamalgam}), may be an equivalence in some cases.

\begin{thm}\label{theorem2.9}
There is a constant $C>0$ such that $\left\|\cdot\right\|_{\alpha}\leq C\left\|\cdot\right\|_{q,p,\alpha}$ whenever $q=\alpha$ or $\alpha=p$.
\end{thm}

In the case $q<\alpha<p$, the space $(L^{q},L^{p})^{\alpha}(X)$ contains properly $L^{\alpha}(X)$ as it appears in the following theorem.

\begin{thm}\label{theorem2.10}
Assume that $1\leq q<\alpha<p\leq\infty$. Then there is a constant $C$ such that
\begin{equation*}
\left\|\cdot\right\|_{q,p,\alpha}\leq C\left\|\cdot\right\|_{\alpha,p}.
\end{equation*}
\end{thm}

The previous result may be strengthened in some cases.

\begin{thm}\label{theorem2.11}
Assume that $1\leq q<\alpha<p$ and that there exists a non decreasing function $\varphi$ on $\left[0,\infty\right)$ and two constants $0<\mathfrak a\leq\mathfrak b<\infty$ such that
\begin{equation}
\mathfrak a\varphi(r)\leq\mu(B_{(x,r)})\leq\mathfrak b \varphi(r),\ x\in X,r>0.\label{normale}
\end{equation}
Then there is $C>0$ such that 
\begin{equation*}
\left\|\cdot\right\|_{q,p,\alpha}\leq C\left\|\cdot\right\|^{\ast}_{\alpha,\infty}.
\end{equation*}
\end{thm}

\medskip

From the doubling condition (\ref{0.05}) and the reverse doubling condition (\ref{revd}), we obtain that the function $\varphi$ appearing in hypothesis (\ref{normale}) satisfies
\begin{equation}
\mathfrak a_{0}r^{D_{\mu}}\leq\varphi(r)\leq\mathfrak b_{0}r^{\delta_{\mu}},\ \ \ r\leq 1,\label{Ip}
\end{equation}
\begin{equation}
\mathfrak a_{1}r^{\delta_{\mu}}\leq\varphi(r)\leq\mathfrak b_{1}r^{D_{\mu}},\ \ \ 1\leq r,\label{Ig}
\end{equation}
where $\mathfrak a_{0},\mathfrak b_{0},\mathfrak a_{1}$ and $\mathfrak b_{1}$ are positive constants.

Notice that Hypothesis (\ref{normale}) is fulfilled for example in the case where $X$ is an Ahlfors $n$ regular metric space, i.e., there is a positive integer $n$ and a constant $C>0$ such that 
$$C^{-1}r^{n}\leq\mu\left(B_{(x,r)}\right)\leq C r^{n},\ x\in X,\ r>0,$$ 
and also in the case where $X$ is a Lie group with polynomial growth equipped with a left  Haar measure $\mu$ and the Carnot-Carath\'eodory metric $d$ associated with a H\"ormander system of left invariant vector fields (see \cite{HMY}, \cite{Ma} and \cite{Va}).

\medskip

The next result shows that the inclusion of $L^{\alpha,\infty}(X)$ into $\left(L^{q},L^{p}\right)^{\alpha}(X)$ is proper.
\begin{thm}\label{theorem2.12}
Under the hypothesis of Theorem \ref{theorem2.11}, we have $(L^{q},L^{p})^{\alpha}(X)\setminus L^{\alpha,\infty}(X)\neq\emptyset$ 
\end{thm}

\section{Auxiliary results}
In order to establish various inclusions between the function spaces we study, we need the following "dyadic cube decomposition" of $X$, proved in \cite{sw}.
\begin{lem}\label{dyadiqueW}
There is $\rho >1$, depending only on $\kappa $ in $%
\left( \ref{0.001}\right) $ (we may take $\rho =8\kappa ^{5}$) ,
such that, given any integer $m$, there exists a family $\left\{(x^{k}_{j},E^{k}_{j}):k\in \mathbb Z,k\geq m,1\leq j<N_{k}\right\}$ where $x^{k}_{j}$ are points of $X$ and $E^{k}_{j}$ subsets of $X$ satisfying:
\begin{enumerate}
\item [(i)]$N_{k}\in\mathbb N^{\ast}\cup\left\{\infty\right\},\ k\geq m$,
\item[(ii)] $B_{\left( x^{k}_{j},\rho ^{k}\right) }\subset
E^{k}_{j}\subset B_{\left( x^{k}_{j},\rho ^{k+1}\right) },\ k\geq m,\ 1\leq j<N_{k},$ 
\item [(iii)] $X=\cup^{N_{k}}_{j=1} E^{k}_{j},$ and $%
E^{k}_{i}\cap E^{k}_{j}=\emptyset $ if $i\neq j$, $k\geq m $,
\item[(iv)] $E^{k}_{j}\subset E^{\ell}_{i}\ $or $E^{k}_{j}\cap E^{\ell}_{i}=\emptyset $, $\ell>k\geq m$, $1\leq j<N_{k}$, $1\leq i<N_{\ell}$.
\end{enumerate}
\end{lem}

The $E^{k}_{j}$  are referred to as dyadic cubes of generation $k$.

\medskip

\begin{notation}\label{4.2}
Given an integer $k\geq m$ and $r>0$, we set
\begin{enumerate}
\item [(i)]$T^{k}_{r}(x)=\left\{i:1\leq i<N_{k}\text{ and } E^{k}_{i}\cap B_{(x,r)}\neq\emptyset\right\},\ x\in X$,
\item[(ii)] $S^{k}_{r}(j)=\left\{i:1\leq i<N_{k}\text{ and } E^{k}_{i}\cap B_{(y,r)}\neq\emptyset\text{ for some }y\in E^{k}_{j}\right\}, 1\leq j<N_{k}$.
\end{enumerate}
\end{notation}

\medskip

Remark that $i\in S^{k}_{r}(j)$ if and only if $j\in S^{k}_{r}(i)$. Inequality (\ref{0.05}) provides us with the following useful estimates on the cardinals $\#(S^{k}_{r}(x))$ and $\#(T^{k}_{r}(x))$ of the sets $S^{k}_{r}(x)$ and $T^{k}_{r}(x)$ respectively.

\medskip

\begin{lem}\label{taille}
Given integers $k\geq m$, $1\leq j<N_{k}$ and $r>0$, we have
\begin{equation}
\mu\left(B_{(y,r)}\right)\leq\mathfrak N_{1}(k,r)\mu\left(E^{k}_{j}\right),\ y\in E^{k}_{j},\label{12}
\end{equation}
\begin{equation}
\mu\left(E^{k}_{i}\right)\leq\mathfrak N_{2}(k,r)\mu\left(E^{k}_{j}\right)\text{ and } \mu\left(E^{k}_{j}\right)\leq\mathfrak N_{2}(k,r)\mu\left(E^{k}_{i}\right),\  i\in S^{k}_{r}(j),\label{13}
\end{equation}
\begin{equation}
\#(T^{k}_{r}(x))\leq\mathfrak N_{2}(k,r),\ x\in X,\label{14}
\end{equation}
and
\begin{equation}
\#(S^{k}_{r}(x))\leq\mathfrak N_{3}(k,r),\ x\in X,\label{15}
\end{equation}
where $\mathfrak N_{1}(k,r)=C_{\mu}\left[\kappa\left(\rho+\frac{r}{\rho^{k}}\right)\right]^{D_{\mu}}$, 
$\mathfrak N_{2}(k,r)=C_{\mu}\left[\kappa\left(2\kappa\rho+\frac{r}{\rho^{k}}\right)\right]^{D_{\mu}}$\\ and 
$\mathfrak N_{3}(k,r)=C_{\mu}\left[\kappa\left((2\kappa^{2}+1)\rho+\frac{r}{\rho^{k}}\right)\right]^{D_{\mu}}\mathfrak N_{2}(k,r)$.
\end{lem}

\proof
\begin{enumerate}
\item[(a)] Inequalities (\ref{12}) and (\ref{13}) are obtained immediately from inequality (\ref{0.05}), the following inclusions:\\
$\bullet\;B_{(x^{k}_{j},\rho^{k})}\subset B_{(x^{k}_{j},\kappa(\rho^{k+1}+r))}$ and $B_{(y,r)}\subset B_{(x^{k}_{j},\kappa(\rho^{k+1}+r))},\ y\in E^{k}_{j}$\\
$\bullet\;E^{k}_{i}\subset B_{(y,\kappa(2\kappa\rho^{k+1}+r))}$ and $B_{(x^{k}_{j},\rho^{k})}\subset B_{(y,\kappa(2\kappa\rho^{k+1}+r))},\ y\in E^{k}_{j}\text{ and } E^{k}_{i}\cap B_{(y,r)}\neq\emptyset,$\\
and the remark stated after Notation \ref{4.2}.
\item[(b)] Lemma \ref{dyadiqueW} (iii) asserts that the $E^{k}_{i}$ $(1\leq i<N_{k})$ are pairwise disjoints. Furthermore we have the following inclusions.\\
$\bullet B_{(x^{k}_{i},\rho^{k})}\subset E^{k}_{i}\subset B_{(x,\kappa(2\kappa\rho^{k+1}+r))}\ \ x\in X\text{ and } i\in T^{k}_{r}(x)$,\\
$\bullet E^{k}_{i}\subset B_{(x^{k}_{j},\kappa[(2\kappa^{2}+1)\rho^{k+1}+r])}\text{ and }  B_{(x^{k}_{j},\rho^{k})}\subset B_{(x^{k}_{j},\kappa[(2\kappa^{2}+1)\rho^{k+1}+r])},\ i\in T^{k}_{r}(j).$

Thus by Inequality (\ref{0.05}), we obtain for $x\in X$
\begin{eqnarray*}
\#(T^{k}_{r}(x))C^{-1}_{\mu}\left[\kappa(2\kappa\rho+\frac{r}{\rho^{k}})\right]^{-D_{\mu}}\mu\left(B_{(x,\kappa(2\kappa\rho^{k+1}+r))}\right)&\leq&\sum_{i\in T^{k}_{r}(x)}\mu\left(B_{(x^{k}_{j},\rho^{k})}\right)\\
&\leq&\mu\left(B_{(x,\kappa(2\kappa\rho^{k+1}+r))}\right)
\end{eqnarray*}
and similarly
\begin{eqnarray*}
\#(S^{k}_{r}(j))\mathfrak N_{2}^{-1}\mu\left(E^{k}_{j}\right)&\leq&\sum_{i\in S^{k}_{r}(j)}\mu\left(E^{k}_{i}\right)\leq\mu\left(B_{(x^{k}_{j},\kappa[(2\kappa^{2}+1)\rho^{k+1}+r])}\right)\\
&\leq& C_{\mu}\left[\kappa((2\kappa^{2}+1)\rho+\frac{r}{\rho^{k}})\right]^{D_{\mu}}\mu\left(E^{k}_{j}\right).
\end{eqnarray*}
Inequalities (\ref{14}) and (\ref{15}) follow.
\end{enumerate}
\epf

\medskip

\begin{lem}\label{lemme3.4}
Assume that $1\leq q,p\leq\infty$, with $p\neq\infty,\;0\leq s$, $m$ and $k$ are integers satisfying $k\geq m$, $1\leq j<N_{k}$ and $2\kappa\rho^{k+1}\leq r$. Then, for any $\mu$-measurable function $f$ on $X$, we have
\begin{equation*}
\mu\left(E^{k}_{j}\right)^{-s}\left\|f\chi_{E^{k}_{j}}\right\|^{p}_{q}\leq\mathfrak N_{1}(k,r)^{s+1}\int_{E^{k}_{j}}\mu\left(B_{(y,r)}\right)^{-s-1}\left\|f\chi_{B_{(y,r)}}\right\|^{p}_{q}d\mu(y)
\end{equation*}
where $\mathfrak N_{1}(k,r)$ is as in Inequality (\ref{12}).
\end{lem}

\proof
Notice that
\begin{equation*}
\underset{E^{k}_{j}}{\inf \text{ess}}\left\|f\chi_{B_{(y,r)}}\right\|^{p}_{q}\leq\mu \left( E^{k}_{j}\right) ^{-1}\int_{E^{k}_{j}}\left\| f\chi_{B_{(y,r)}}\right\|^{p}_{q}d\mu(y) 
\end{equation*}
with equality only when $\left\|f\chi_{B_{(y,r)}}\right\|_{q}$ is a constant almost everywhere on $E^{k}_{j}$. 
Thus, there is an element $y^{k}_{j}$ of $E^{k}_{j}$ such that
 \begin{equation*}
 \left\|f\chi_{B_{(y^{k}_{j},r)}}\right\|^{p}_{q}\leq\mu \left( E^{k}_{j}\right) ^{-1}\int_{E^{k}_{j}}\left\| f\chi_{B_{(y,r)}}\right\|^{p}_{q}d\mu(y).
\end{equation*}
Since $E^{k}_{j}$ is included in $B_{(y,r)}$ for every $y$ in $E^{k}_{j}$, we have 
\begin{equation*}
\mu\left(E^{k}_{j}\right)^{-s}\left\|f\chi_{E^{k}_{j}}\right\|^{p}_{q}\leq\mu\left(E^{k}_{j}\right)^{-s}\left\|f\chi_{B_{(y^{k}_{j},r)}}\right\|^{p}_{q}\leq\mu\left(E^{k}_{j}\right)^{-s-1}\int_{E^{k}_{j}}\left\|f\chi_{B_{(y,r)}}\right\|^{p}_{q}d\mu(y).
\end{equation*}
The result follows from inequality (\ref{12}).
\epf

\medskip

We shall use the following result which may be viewed as a generalization of the Young inequality in a space without group structure.

\medskip

\begin{lem}
\label{younginequality} Let $\beta ,t$ and $\gamma$ be elements of $%
\left[ 1\ ,\infty \right] $ such that $\frac{1}{\gamma}=\frac{1}{\beta}+\frac{1}{t}-1$ and $K\left( x,y\right)$ a positive
kernel on $X$. There is a constant $C>0$ such that 
\begin{equation*}
\left\|Tg\right\|_{\gamma}\leq C\left\| \left\|
K\right\| _{\beta }\right\| _{\infty }\left\| g\right\|^{\ast} _{t,\gamma },\ g\in L_{0}(X),
\end{equation*}
where 
\begin{equation*}
Tg(y)=\int_{X}g(x)K(x,y)d\mu(x),
\end{equation*}
 and 
 \begin{equation*}
 \left\| \left\| K\right\| _{\beta }\right\| _{\infty
}=\max \left( \underset{y\in X}{\sup ess}\left\| K\left( .,y\right) \right\|
_{\beta };\underset{x\in X}{\sup ess}\left\| K\left( x,.\right) \right\|
_{\beta }\right). 
\end{equation*}
\end{lem}

\proof
\begin{enumerate}
\item [1)]Let $g\in L_{0}\left( X\right) $ and put $\widetilde{g}(y)=\int_{X}\left| g(x)\right| K\left( x,y\right) d\mu\left(x\right) $.
We claim that 
\begin{equation*}
\left\|\tilde{g}\right\| ^{\ast}_{\gamma,\infty }\leq C\left\| \left\|
K\right\| _{\beta }\right\| _{\infty }\left\| g\right\| _{t,\infty }^{\ast
}.
\end{equation*}
If $g\notin L^{t,\infty }\left( X \right) $ or $\left\| g\right\|^{\ast}
_{t,\infty }=0,$ or $\left\| \left\| K\right\| _{\beta }\right\|
_{\infty }\in\left\{0,\infty\right\}$ then the claim is trivially verified. So we assume that $0<\left\| g\right\|^{\ast} _{t,\infty }<\infty $ and $%
0<\left\| \left\| K\right\| _{\beta }\right\| _{\infty }<\infty $. Define
\begin{equation*}
g_{_{1}}(x)=\left\{ 
\begin{array}{ll}
g(x) & \text{if \ }\left| g(x)\right| \leq M \\ 
0 & \text{if not}%
\end{array}
\right. \text{ and } g_{_{2}}(x)=g(x)-g_{1}(x),\ x\in X,
\end{equation*}
 where $M$ is a positive real number to be specified later. For $\alpha>0$, we have  $\lambda _{\widetilde{g}}\left( \alpha \right) \leq \lambda _{%
\widetilde{g}_{_{1}}}\left( \frac{\alpha }{2}\right) +\lambda _{\widetilde{g%
}_{_{2}}}\left( \frac{\alpha }{2}\right)$ since $\widetilde{g}\leq \widetilde{g}_{_{1}}+\widetilde{g}_{_{2}}$.

\begin{enumerate}
\item[$a)$] We can estimate $\lambda _{\widetilde{g}_{_{1}}}\left( \frac{\alpha }{2}\right)$ as follows:
\begin{eqnarray*}
\int_{X}\left| g_{_{1}}\left( x\right) \right| ^{\beta'}d\mu (x) &=&\beta'\int_{0}^{\infty }s^{\beta'-1}\lambda
_{g_{_{1}}}\left( s\right) ds\leq\beta'\int_{0}^{M}s^{\beta'-1}\lambda _{g}\left( s\right) ds \\ 
& \leq& \beta'\left( \int_{0}^{M}s^{\beta'-1-t}ds\right)
\left( \left\| g\right\| _{t,\infty }^{\ast }\right) ^{t}=\frac{\beta'}{\beta'-t}M^{\beta'-t}\left(\left\|g\right\|^{\ast}_{t,\infty}\right)^{t}.
\end{eqnarray*}
So, 
\begin{eqnarray*}
\left| \widetilde{g}_{_{1}}\left( y\right) \right| &=& \int_{X}\left|
g_{_{1}}(x)\right| K(x,y)d\mu \left( x\right) \leq \left( \int_{X}\left|
g_{_{1}}(x)\right| ^{\beta ^{\prime }}d\mu \left( x\right) \right) ^{\frac{1%
}{\beta ^{\prime }}}\left( \int_{X}K^{\beta }(x,y)d\mu \left( x\right)
\right) ^{\frac{1}{\beta }} \\ 
&\leq& \left(\frac{\gamma}{t}\right)^{\frac{1}{\beta'}}M^{\frac{t}{\gamma}}\left(\left\|g\right\|^{\ast}_{t,\infty}\right)^{-\frac{t}{\beta'}} \left\| \left\| K\right\| _{\beta }\right\|
_{\infty }.
\end{eqnarray*}
Let us choose  
\begin{equation*}
M=\left( \frac{\alpha }{2}\right) ^{\frac{\gamma}{t}}\left( \frac{t}{\gamma}\right)
^{\frac{\gamma}{t\beta'}}\left( \left\| g\right\|^{\ast
} _{t,\infty }\right) ^{-\frac{\gamma}{\beta'}}\left\| \left\| K\right\| _{\beta
}\right\|^{-\frac{\gamma}{t}} _{\infty }.
\end{equation*}
We have $\left\|\tilde{g}_{1}\right\|_{\infty}\leq\frac{\alpha}{2}$ and therefore $\lambda _{\widetilde{g}_{_{1}}}\left( \frac{\alpha }{2}%
\right) =0.$

\item[$b)$]  We also have the following estimate of $\lambda _{\widetilde{g}_{_{2}}}\left( \frac{\alpha }{2}\right)$:
\begin{eqnarray*}
\int_{X}\left| g_{_{2}}(x)\right| d\mu \left( x\right) & =&\int_{0}^{\infty
}\lambda _{g_{_{2}}}(s)ds\leq \int_{0}^{M}\lambda
_{g}(M)ds+\int_{M}^{\infty }\lambda _{g}(s)ds \\ 
& \leq& M^{1-t}\left( \left\| g\right\| _{t,\infty }^{\ast }\right)
^{t}+\left( \int_{M}^{\infty }s^{-t}ds\right) \left( \left\| g\right\|
_{t,\infty }^{\ast }\right) ^{t}=\left( \frac{t}{t-1}\right) M^{1-t}\left(
\left\| g\right\| _{t,\infty }^{\ast }\right) ^{t}.
\end{eqnarray*}
Therefore, 
\begin{eqnarray*}
\lambda _{\widetilde{g}_{_{2}}}\left( \frac{\alpha }{2}\right)& \leq& \left( \frac{2}{\alpha }\right) ^{\beta }\int_{\left\{ u\in X:\left| 
\widetilde{g}_{_{2}}(u)\right| >\frac{\alpha }{2}\right\} }\left(
\int_{X}\left| g_{_{2}}(x)\right| K(x,y)d\mu (x)\right) ^{\beta }d\mu \left(
y\right) \\ 
& \leq& \left( \frac{2}{\alpha }\right) ^{\beta }\left[ \int_{X}\left|
g_{_{2}}(x)\right| \left( \int_{\left\{ u\in X:\left| \widetilde{g}%
_{_{2}}(u)\right| >\frac{\alpha }{2}\right\} }K^{\beta }(x,y)d\mu (y)\right)
^{\frac{1}{\beta }}d\mu \left( x\right) \right] ^{\beta } \\ 
& \leq &\left( \frac{2}{\alpha }\right) ^{\beta }\left\| \left\| K\right\|
_{\beta }\right\| _{\infty }^{\beta }\left[ \int_{X}\left|
g_{_{2}}(x)\right| d\mu \left( x\right) \right] ^{\beta } \\ 
& \leq &\left( \frac{2}{\alpha }\right) ^{\beta }\left\| \left\| K\right\|
_{\beta }\right\| _{\infty }^{\beta }\left[ \left( \frac{t}{t-1}\right)
M^{1-t}\left( \left\| g\right\| _{t,\infty }^{\ast }\right) ^{t}\right]
^{\beta }\leq \left(C\alpha ^{-1}\left\| \left\| K\right\| _{\beta }\right\|
_{\infty }\left\| g\right\| _{t,\infty }^{\ast }\right) ^{\gamma},%
\end{eqnarray*}
with $C=\left( 2\right) ^{\gamma}\left(\frac{t}{t-1}\right) ^{\beta }\left( 
\frac{t}{\gamma}\right) ^{\frac{\gamma}{t\beta'}\left( 1-t\right) \beta }$.
\end{enumerate}

From a) and b) we get
\begin{equation*}
\lambda _{\widetilde{g}}\left( \alpha
\right) \leq \left(C\alpha ^{-1}\left\| \left\| K\right\| _{\beta }\right\|
_{\infty } \left\| g\right\| _{t,\infty }^{\ast }\right) ^{\gamma}.
\end{equation*}
As this inequality is true for $\alpha>0$, we have
\begin{equation*}
\left\|Tg\right\|^{\ast}_{\gamma,\infty}\leq C\left\|\left\|K\right\|_{\beta}\right\|_{\infty}\left\|g\right\|^{\ast}_{t,\infty}.
\end{equation*}
\item [2)] Notice that $T$ is a linear operator. Therefore, the result follows from $1)$ and Stein interpolation theorem (see \cite{SWe}).
\end{enumerate}
\epf

\section{Proof of the main results}
Throughout this paragraph, for every $r>0$, $m_{r}$ denotes the unique integer 
 which verifies
\begin{equation}
\rho ^{m_{r}+1}\leq \frac{r}{2\kappa }<\rho ^{m_{r}+2}.  \label{0.002}
\end{equation}
Notice that the constants in Lemma \ref{taille} satisfy
\begin{equation}
\mathfrak N_{1}(m_{r},r)\leq C_{\mu}\left[\kappa\rho(1+2\kappa\rho)\right]^{D_{\mu}}=\mathfrak N_{1},\label{13'}
\end{equation} 
\begin{equation}
\mathfrak N_{2}(m_{r},r)\leq C_{\mu}\left[2\kappa^{2}\rho(1+\rho)\right]^{D_{\mu}}=\mathfrak N_{2},\label{14'}
\end{equation} 
and
\begin{equation}
\mathfrak N_{3}(m_{r},r)\leq C_{\mu}\left[\kappa\rho(2\kappa^{2}+2\kappa\rho+1)\right]^{D_{\mu}}\mathfrak N_{2}=\mathfrak N_{3}.\label{15'}
\end{equation} 

\medskip

\proof [\bf{Proof of Theorem \ref{theorem2.2}}]
 Let $f$ be a $\mu$-measurable function on $X$ such that $
_{r}\left\| f\right\| _{q,p,\alpha}=0.$ Since balls in $X$ have positive measure, $\left\|\mu(B_{(\cdot,r)})^{\frac{1}{\alpha}-\frac{1}{p}-\frac{1}{q}} \left\| f\chi _{_{B\left(\cdot,r\right) }}\right\| _{q}\right\|_{p}=0$ implies that there exists a $\mu$-null subset $E$ of $X$ such that
\begin{equation*}
\left\| f\chi _{_{B\left(\cdot,r\right) }}\right\| _{q}=0\text{ in }X\setminus E.
\end{equation*}
Similarly, for any $y$ in $X\setminus E$, there exists a $\mu$-null subset $F_{y}$ of $X$ out of which $f\chi_{B_{(y,r)}}=0$.
For $1\leq j<N_{m_{r}}$, the intersection of $X\setminus E$ and $B_{\left(x^{m_{r}}_{j},\rho ^{m_{r}+1}\right) }$ is non void. So we may pick in it an element $y_{j}$. Since $E^{m_{r}}_{j}\subset B_{\left(x_{j}^{m_{r}},\rho ^{m_{r}+1}\right) }\subset B_{\left( y_{j},r\right) }$, we have $X=\underset{j=1}{\overset{N_{m_{r}}}{\cup }}E^{m_{r}}_{j}=\underset{j=1}{\overset{N_{m_{r}}}{\cup }}B_{\left( y_{j},r\right) }.$ Setting $F=\cup^{N_{m_{r}}}_{j=1}F_{y_{j}}$, we have $f=0$ in $X\setminus F$. The result follows from the fact that $\mu(F)=0$.
\epf

\medskip

\proof[\bf{Proof of Theorem \ref{eamalgamX}}]
It is clear from Theorem \ref{theorem2.2} and the definition of $_{r}\left\| \cdot\right\| _{q,p,\alpha}$ that $(L^{q},L^{p})^{\alpha}_{r}(X)$ is a complex vector space and $_{r}\left\| \cdot\right\| _{q,p,\alpha}$ is a norm on it. All we need to prove is completeness.

Let $\left( f_{n}\right) _{n>0}$ be a sequence of elements of $\left(
L^{q},L^{p}\right)^{\alpha} _{r}\left( X \right) $ such that $\underset{n>0}{%
\sum }\ _{r}\left\| f_{n}\right\| _{q,p,\alpha}<\infty .$

Since $\underset{n>0}%
{\sum }\ \left\| \mu\left(B_{(\cdot,r)}\right)^{\frac{1}{\alpha}-\frac{1}{p}-\frac{1}{q}}\left\| f_{n}\chi _{_{B_{\left( \cdot,r\right) }}}\right\|
_{q}\right\| _{p}=\underset{n>0}{\sum }\ _{r}\left\| f_{n}\right\| _{q,p,\alpha}<\infty,$ there exists a $\mu -$null  subset $E$ of $X$ out of which 
\begin{equation*}
\sum _{n>0} \left\| f_{n}\chi _{_{B_{\left( y,r\right)
}}}\right\| _{q}<\infty.
\end{equation*}
Therefore, for any element $y$ of $X\setminus E$, there is a $\mu -$null subset $F_{y}$ of $X$ out of which $\underset{n>0}{\sum }f_{n} \chi _{_{B_{\left( y,r\right)}}}$ converges absolutely. Arguing as in the proof of Theorem \ref{theorem2.2}, we shall obtain a $\mu -$null subset $F$ of $X$ such that $\sum_{n>0}f_{n}$ converges absolutely on $X\setminus F$.
Define 
\begin{equation*}
f(x)=\left\{ 
\begin{array}{ll}
\underset{n>0}{\sum }f_{n}\left( x\right) & \text{if }x\in X\setminus F \\ 
0 & \text{otherwise}%
\end{array}
\right. .
\end{equation*}
We have 
\begin{equation*}
_{r}\left\| f\right\| _{q,p,\alpha}\leq \underset{n>0}{\sum }\ _{r}\left\|
f_{n}\right\| _{q,p,\alpha}<\infty .
\end{equation*}
In addition, for any positive integer $n$ and any element $y$ of $X$, 
\begin{equation*}
\left\| f\chi _{_{B_{\left( y,r\right) }}}-\underset{k=1}{\overset{n}{\sum }}%
f_{k}\chi _{_{B_{\left( y,r\right) }}}\right\| _{q}\leq \underset{k>n}{\sum }%
\left\| f_{k}\chi _{_{B_{\left( y,r\right) }}}\right\| _{q}.
\end{equation*}
Therefore 
\begin{equation*}
_{r}\left\|f-\underset{k=1}{\overset{n}{\sum}}f_{k}\right\| _{q,p,\alpha}\leq 
\underset{k>n}{\sum }\ _{r}\left\| f_{k}\right\| _{q,p,\alpha}.
\end{equation*}

Thus
 $\underset{n>0}{\sum }f_{n}$ converges to $f$ in $\left( L^{q},L^{p}\right)^{\alpha}
_{r}\left( X\right) $.
\epf

The norm $_{r}\left\|\cdot\right\| _{q,p,\alpha}$ is not easy to be used. The following proposition provides us with an equivalent norm.

\begin{prop}\label{normequivalente}
Let $f$ be any $\mu$-measurable function on $X$, and $r>0$. Put
\begin{equation*}
\left\| f\right\| _{q,p,\alpha}^{dm_{r}}=\left\{ 
\begin{array}{lll}
\left[\sum^{N_{m_{r}}}_{j=1}\left(\mu\left(E^{m_{r}}_{j}\right)^{\frac{1}{\alpha}-\frac{1}{q}}\left\| f\chi_{_{E_{j}^{m_{r}}}}\right\| _{q}\right)^{p}\right]^{\frac{1}{p}} & \text{ if } & p<\infty \\ 
\underset{1\leq j<N_{m_{r}}}{\sup}\mu\left(E^{m_{r}}_{j}\right)^{\frac{1}{\alpha}-\frac{1}{q}}\left\| f\chi _{_{E_{j}^{m_{r}}}}\right\|
_{q} & \text{if} & p=\infty%
\end{array}
\right. .
\end{equation*}
Then, there are positive constants $C_{1}$and $C_{2}$, not depending on $f$ and $r$, such that  
\begin{equation}
C_{1}\;_{r}\left\| f\right\| _{q,p,\alpha}\leq \left\|
f\right\| _{q,p,\alpha}^{dm_{r}}\leq C_{2}\;_{r}\left\|
f\right\| _{q,p,\alpha}.\label{equinorm}
\end{equation}
\end{prop}

\proof
Let $f$ be any $\mu$-measurable function on $X$ and $r>0$. 
\begin{enumerate}
\item[{\it 1$^{st}$ case. }] We suppose that $p<\infty .$
\begin{enumerate}
\item [a)] We have
\begin{eqnarray*}
_{r}\left\| f\right\| _{q,p,\alpha}^{p} & =&\int_{X}\left\{ \mu\left(B_{(y,r)}\right)^{\frac{q}{\alpha}-\frac{q}{p}-1}\int_{X}\left( \left|
f\right| ^{q}\chi _{_{B_{\left( y,r\right) }}}\right) (x)d\mu (x)\right\}^{\frac{p}{q%
}}d\mu \left( y\right)  \\ 
& =&\underset{j=1}{\overset{N_{m_{r}}}{\sum }}\int_{E_{j}^{m_{r}}}\left\{ 
\underset{i\in T^{m_{r}}_{r}(y)}{\sum }\mu\left(B_{(y,r)}\right)^{\frac{q}{\alpha}-\frac{q}{p}-1}\int_{E^{m_{r}}_{i}}\left(
\left| f\right| ^{q}\chi _{_{B_{\left( y,r\right) }}}\right) (x)d\mu
(x)\right\} ^{\frac{p}{q}}d\mu \left( y\right) \\ 
& \leq& \mathfrak N^{\frac{p}{q}-1}_{2}\underset{j=1}{\overset{N_{m_{r}}}{\sum }}%
\int_{E_{j}^{m_{r}}}\mu\left(B_{(y,r)}\right)^{\frac{p}{\alpha}-\frac{p}{q}-1}\underset{i\in T^{m_{r}}_{r}(y)}{\sum }\left[
\int_{E^{m_{r}}_{i}}\left( \left| f\right| ^{q}\chi _{_{B_{\left(
y,r\right) }}}\right) (x)d\mu (x)\right] ^{\frac{p}{q}}d\mu \left( y\right),
\end{eqnarray*}
according to Inequalities (\ref{14}) and (\ref{14'}). As $2\kappa\rho^{m_{r}+1}\leq r,$ we have $E^{m_{r}}_{i}\subset B_{(y,2\kappa r)}$ for $i\in T^{m_{r}}_{r}(y)$ and therefore by Inequality (\ref{0.05}), 
\begin{equation}
\mu\left(E^{m_{r}}_{i}\right)\leq C_{\mu}(2\kappa)^{D_{\mu}}\mu\left(B_{(y,r)}\right), \ i\in T^{m_{r}}_{r}(y).\label{23} 
\end{equation}
Taking into account Inequalities (\ref{23}), (\ref{13}) and (\ref{14'}), we obtain 
\begin{equation*} 
_{r}\left\| f\right\| _{q,p,\alpha}^{p}\leq C\sum^{N_{m_{r}}}_{j=1}\sum_{i\in S^{m_{r}}_{r}(j)}\mu(E^{m_{r}}_{i})^{\frac{p}{\alpha}-\frac{p}{q}}\left\| f\chi _{_{E^{m_{r}}_{i}}}\right\| _{q}^{p}.
\end{equation*}
 So by Inequalities (\ref{15}) and (\ref{15'}) we get
\begin{equation*} 
_{r}\left\| f\right\| _{q,p,\alpha}^{p}\leq C\mathfrak N_{3}\sum^{N_{m_{r}}}_{i=1}\mu(E^{m_{r}}_{i})^{\frac{p}{\alpha}-\frac{p}{q}}\left\| f\chi _{_{E^{m_{r}}_{i}}}\right\| _{q}^{p}\leq C\mathfrak N_{3}\left( \left\| f\right\|
_{q,p,\alpha}^{d_{m_{r}}}\right) ^{p}.%
\end{equation*}
\item[b)]
Notice that if $_{r}\left\| f\right\| _{q,p,\alpha}=\infty$, then (\ref{equinorm}) follows trivially from the above inequality. Let us assume that $_{r}\left\| f\right\| _{q,p,\alpha}<\infty$. For $1\leq j<N_{m_{r}}$, we have   
\begin{equation*}
\mu\left(E^{m_{r}}_{j}\right)^{\frac{1}{\alpha}-\frac{1}{q}}\left\|f\chi_{E^{m_{r}}_{j}}\right\|^{p}_{q}\leq\mathfrak N_{1}^{\frac{1}{q}-\frac{1}{\alpha}+1}\int_{E^{m_{r}}_{j}}\mu\left(B_{(y,r)}\right)^{\frac{1}{\alpha}-\frac{1}{q}-1}\left\|f\chi_{B_{(y,r)}}\right\|^{p}_{q}d\mu(y),
\end{equation*}
according to Lemma \ref{lemme3.4}. As the $E^{m_{r}}_{j}$ $(1\leq j<N_{m_{r}})$ are pairwise disjoints, this implies
\begin{equation*}
\left\|f\right\|^{d_{m_{r}}}_{q,p,\alpha}\leq\mathfrak N^{\frac{1}{p}(\frac{1}{q}-\frac{1}{\alpha}+1}_{1}\ _{r}\left\|f\right\|_{q,p,\alpha}.
\end{equation*}
\end{enumerate}
\item[{\it 2$^{nd}$case.}] We suppose that $p=\infty $.
\begin{enumerate}
\item  We have 
\begin{eqnarray*}
_{r}\left\| f\right\| _{q,\infty ,\alpha} & =&\underset{y\in X}{\sup }
\left[ \underset{j\in T^{m_{r}}_{r}(y)}{\sum }\mu\left(B_{(y,r)}\right)^{\frac{q}{\alpha}-1}
\int_{E_{j}^{m_{r}}}\left| f(x)\chi _{_{B_{\left( y,r\right) }}}\left(
x\right) \right| ^{q}d\mu (x)\right] ^{\frac{1}{q}} \\ 
& \leq &\left[C_{\mu}(2\kappa)^{D_{\mu}}\right]^{\frac{1}{\alpha}-\frac{1}{q}}\underset{y\in X}{\sup }
 \underset{j\in T^{m_{r}}_{r}(y)}{\sum }\mu\left(E^{m_{r}}_{j}\right)^{\frac{1}{\alpha}-\frac{1}{q}}
\left\| f\chi _{_{E^{m_{r}}_{j}}}\right\|_{q}\\
&\leq&\left[C_{\mu}(2\kappa)^{D_{\mu}}\right]^{\frac{1}{\alpha}-\frac{1}{q}}\mathfrak N_{2}\left\| f\right\| _{q,\infty,\alpha }^{d_{m_{r}}},
\end{eqnarray*}
according to Inequalities (\ref{23}), (\ref{14}) and (\ref{14'}).
\item From Inequalities (\ref{12}) and (\ref{13'}) we have 
\begin{equation*}
\mu(E^{k}_{j})^{\frac{1}{\alpha}-\frac{1}{q}}\leq\mathfrak N_{1}\mu(B_{(y,r)})^{\frac{1}{\alpha}-\frac{1}{q}},\ \ \ 1\leq j<N_{m_{r}}\text{ and } y\in E^{m_{r}}_{j}
\end{equation*}
and therefore  
\begin{equation*}
\left\| f\right\| _{q,\infty,\alpha }^{d_{m_{r}}}  \leq\mathfrak N_{1}\underset{1\leq
j<N_{m_{r}}}{\sup}\sup_{y\in E^{m_{r}}_{j}}\mu(B_{(y,r)})^{\frac{1}{\alpha}-\frac{1}{q}}\left\|f\chi_{E^{m_{r}}_{j}}\right\|_{q}.
\end{equation*}
As $2\kappa\rho^{m_{r}+1}\leq r$, we have 
\begin{equation*}
E^{m_{r}}_{j}\subset B_{(y,r)},\ \ \ 1\leq j<N_{m_{r}}\text{ and } y\in E^{m_{r}}_{j}.
\end{equation*}
Thus,
\begin{equation*}
\left\| f\right\| _{q,\infty,\alpha }^{d_{m_{r}}}\leq\mathfrak N_{1}\sup_{y\in X}\mu(B_{(y,r)})^{\frac{1}{\alpha}-\frac{1}{q}}\left\|f\chi_{B_{(y,r)}}\right\|_{q}=\mathfrak N_{1}\;_{r}\left\|f\right\|_{q,\infty,\alpha}.
\end{equation*}
\end{enumerate}

\item[{\it 3$^{rd}$case.}] For $q=p=\infty$, it is clear that
\begin{equation}
_{r}\left\|f\right\|_{\infty,\infty,\infty}=\left\|f\right\|_{\infty}=\left\|f\right\|^{dm_{r}}_{\infty,\infty,\infty}.
\end{equation}
\end{enumerate}
\epf

\medskip

\proof[\bf{Proof of Theorem \ref{comparaisonnormamalgame}}] 
\begin{enumerate}
\item[a)] Inequality (\ref{9}) is an immediate consequence of H\"older inequality.
\item[b)] Observe that as $0<p_{1}<p_{2}<\infty$, we have for any sequence $(a_{j})_{1\leq j}$ of nonnegative numbers ,
\begin{equation*}
\sup_{1\leq j}a_{j}\leq \left(\sum^{\infty}_{j=1}a^{p_{2}}_{j}\right)^{\frac{1}{p_{2}}}\leq\left(\sum^{\infty}_{j=1}a^{p_{1}}_{j}\right)^{\frac{1}{p_{1}}}
\end{equation*}
and therefore
\begin{equation*}
\left\|\cdot\right\|^{dm_{r}}_{q,\infty,\alpha}\leq\left\|\cdot\right\|^{dm_{r}}_{q,p_{2},\alpha}\leq\left\|\cdot\right\|^{dm_{r}}_{q,p_{1},\alpha}.
\end{equation*}
\end{enumerate}
Inequality (\ref{10}) follows from these inequalities and Proposition \ref{normequivalente}.
\epf

\medskip

\proof[\bf{Proof of Theorem \ref{comparaisonnormamalgamelebesgue}}]
Let $f$ be any $\mu$-measurable function on $X$.
\begin{enumerate}
\item[{\it 1$^{rst}$ case.}] We suppose that $p=\infty.$

By H\"older inequality we have
\begin{equation*}
\;_{r}\left\|f\right\|_{q,\infty,\alpha}\leq\sup_{y\in X}\left\|f\chi_{B_{(y,r)}}\right\|_{\alpha}\leq\left\|f\right\|_{\alpha}.
\end{equation*}

\item[{\it 2$^{nd}$ case.}] We suppose that $p<\infty.$  Then we have

\begin{equation*}
\left\|f\right\|^{dm_{r}}_{q,p,\alpha}\leq\left[\sum^{N_{m_{r}}}_{j=1}\left(\mu(E^{m_{r}}_{j})^{\frac{1}{\alpha}-\frac{1}{q}}\left\|f\chi_{E^{m_{r}}_{j}}\right\|_{q}\right)^{p}\right]^{\frac{1}{p}}\leq\left(\sum^{N_{m_{r}}}_{j=1}\left\|f\chi_{E^{m_{r}}_{j}}\right\|^{\alpha}_{\alpha}\right)^{\frac{1}{\alpha}}\leq\left\|f\right\|_{\alpha}
\end{equation*}
according to H\"older inequality, the fact that $0<\alpha\leq p<\infty$ and the pairwise disjointness of the $E^{m_{r}}_{j}$ ($1\leq j<N_{m_{r}}$). From this inequality and Proposition \ref{normequivalente} we obtain (\ref{11}).
\end{enumerate}
\epf

\medskip

\proof[\bf{Proof of Theorem \ref{theorem2.9}}]\ 

\begin{enumerate}
\item[{\it 1$^{rst}$ case.}] We suppose that $q=\alpha=p$. 

It is clear from Proposition \ref{normequivalente} that there is a constant $C_{2}$, not depending on $f$, such that
\begin{equation*}
\left\|f\right\|_{\alpha}=\left\|f\right\|^{dm_{r}}_{\alpha,\alpha,\alpha}\leq C_{2}\;_{r}\left\|f\right\|_{\alpha,\alpha,\alpha},\ \ r>0
\end{equation*}
and therefore
\begin{equation*}
\left\|f\right\|_{\alpha}\leq C_{2}\left\|f\right\|_{\alpha,\alpha,\alpha}.
\end{equation*}
\item[{\it 2$^{nd}$ case.}] We suppose that $q=\alpha<p=\infty$.

For any element $y$ of $X$ formula (\ref{namalgamX}) yields
\begin{equation*}
\left\|f\chi_{B_{(y,r)}}\right\|_{\alpha}\leq\;_{r}\left\|f\right\|_{\alpha,\alpha,\infty}\leq\left\|f\right\|_{\alpha,\alpha,\infty},\ r>0
\end{equation*}
and therefore 
\begin{equation*}
\left\|f\right\|_{\alpha}=\lim_{r\rightarrow\infty}\left\|f\chi_{B_{(y,r)}}\right\|_{\alpha}|\leq\left\|f\right\|_{\alpha,\alpha,\infty}.
\end{equation*}
\item[{\it 3$^{rd}$ case.}] We suppose that $q=\alpha<p<\infty$. For $y\in X$ and $r>0$, we have
\begin{eqnarray*}
\left\|f\chi_{B_{(y,r)}}\right\|_{\alpha}&=&\left(\sum^{N_{m_{r}}}_{j=1}\int_{E^{m_{r}}_{j}}\left|f(x)\right|^{\alpha}\chi_{B_{(y,r)}}(x)d\mu(x)\right)^{\frac{1}{\alpha}}\\
&=&\left(\sum_{j\in T^{m_{r}}_{r}(y)}\int_{X}\left|(f\chi_{E^{m_{r}}_{j}})(x)\right|^{\alpha}\chi_{B_{(y,r)}}(x)d\mu(x)\right)^{\frac{1}{\alpha}}\\
&\leq&\left(\sum_{j\in T^{m_{r}}_{r}(y)}\left\|f\chi_{E^{m_{r}}_{j}}\right\|^{\alpha}_{\alpha}\right)^{\frac{1}{\alpha}}\leq\mathfrak N^{\frac{1}{\alpha}-\frac{1}{p}}_{2}\left(\sum_{j\in T^{m_{r}}_{r}(y)}\left\|f\chi_{E^{m_{r}}_{j}}\right\|^{p}_{\alpha}\right)^{\frac{1}{p}}
\end{eqnarray*}
according to Inequalities (\ref{13}) and (\ref{13'}). So by Proposition \ref{normequivalente}, we get
\begin{equation*}
\left\|f\chi_{B_{(y,r)}}\right\|_{\alpha}\leq\mathfrak N^{\frac{1}{\alpha}-\frac{1}{p}}_{2} C_{2}\;_{r}\left\|f\right\|_{\alpha,\alpha,p},\ \ y\in X,\ r>0
\end{equation*}
and therefore
\begin{equation*}
\left\|f\right\|_{\alpha}\leq\mathfrak N^{\frac{1}{\alpha}-\frac{1}{p}}_{2} C_{2}\;_{r}\left\|f\right\|_{\alpha,\alpha,p}.
\end{equation*}
\item[{\it 4$^{th}$ case.}] We suppose that $q<\alpha=p.$
We assume that $\left\| f\right\| _{q,p,p}<\infty ,$ since otherwise the result follows from Theorem \ref{theorem2.8}. For $r>0,$ put  
\begin{equation*}
f_{r}(x)=\mu\left(B_{(y,r)}\right)^{-\frac{1}{q}}\left\|f\chi_{B_{(y,r)}}\right\|_{q}.
\end{equation*}
On one hand, we have for $\mu-$almost every $x$ in $X$,
\begin{equation*}
\left| f(x)\right|=\underset{r\rightarrow 0}{\lim }f_{r}(x) \leq \left\| f\right\| _{q,\infty ,\infty }.
\end{equation*}
 Consequently 
\begin{equation*}
\left\| f\right\| _{\infty }\leq \left\| f\right\|
_{q,\infty ,\infty }.
\end{equation*}
On the other hand, 
\begin{equation*}
\left[ \int_{X}f_{r}^{\text{ }p}(x)d\mu (x)\right] ^{\frac{1}{p}}\leq
C\left\| f\right\| _{q,p,p}.
\end{equation*}
So, according to Fatou's lemma, $\left| f\right| ^{p}$ is integrable and $%
\left\| f\right\| _{p}\leq C\left\| f\right\| _{q,p,p}.$
\end{enumerate}
\epf

\medskip

\proof[\bf{Proof of Theorem \ref{theorem2.10}}]
Let $\frac{1}{\beta}=1-\frac{q}{\alpha}+\frac{q}{p}$, $f$ a $\mu$-measurable function on $X$ and $r>0$.

We have $1<\beta,$ $\frac{\alpha}{q}<\infty$ and $\frac{q}{p}=\frac{1}{\beta}+\frac{q}{\alpha}-1$. Put
\begin{equation*}
K(x,y)=\mu\left(B_{(x,r)}\right)^{-\frac{1}{\beta}}\chi_{B_{(y,r)}}(x),\ \ x,y\in X
\end{equation*}
and
\begin{equation*}
Tg(y)=\int_{X}g(x)K(x,y)d\mu(y),\ \ g\in L_{0}(X).
\end{equation*}
If $x\in B_{(y,r)}$ then $B_{(y,r)}\subset B_{(x,2\kappa r)}$ and therefore $\mu\left(B_{(y,r)}\right)^{-1}\leq C_{\mu}(2\kappa)^{D_{\mu}}\mu\left(B_{(x,r)}\right)^{-1}$. 

Thus 
\begin{equation*}
\left(\int_{X}\left|K(x,y)\right|^{\beta}d\mu(y)\right)^{\frac{1}{\beta}}=\left(\int_{X}\mu\left(B_{(y,r)}\right)^{-1}\chi_{B_{(x,r)}}(y)d\mu(y)\right)^{\frac{1}{\beta}}\leq C_{\mu}(2\kappa)^{D_{\mu}},
\end{equation*}
and
\begin{equation*}
\left(\int_{X}\left|K(x,y)\right|^{\beta}d\mu(x)\right)^{\frac{1}{\beta}}=\left(\int_{X}\mu\left(B_{(y,r)}\right)^{-1}\chi_{B_{(y,r)}}(x)d\mu(x)\right)^{\frac{1}{\beta}}=1.
\end{equation*}
By Lemma \ref{younginequality}, there is a constant $C$ such that 
\begin{equation*}
\left\|T\left(\left|f\right|^{q}\right)\right\|_{\frac{p}{q}}\leq C\left\|\left|f\right|^{q}\right\|^{\ast}_{\frac{\alpha}{q},\frac{p}{q}}.
\end{equation*}
Furthermore 
\begin{equation*}
_{r}\left\|f\right\|_{q,p,\alpha}=\left[\int_{X}\left(T(\left|f\right|^{q})(y)\right)^{\frac{p}{q}}d\mu(y)\right]^{\frac{1}{p}}=\left(\left\|T(\left|f\right|^{q})\right\|_{\frac{p}{q}}\right)^{\frac{1}{q}}.
\end{equation*}
Thus 
\begin{equation*}
_{r}\left\|f\right\|_{q,p,\alpha}\leq\left(C\left\|\left|f\right|^{q}\right\|^{\ast}_{\frac{\alpha}{q},\frac{p}{q}}\right)^{\frac{1}{q}}=C^{\frac{1}{q}}\left\|f\right\|^{\ast}_{\alpha,p}.
\end{equation*}
The result follows.
\epf

\medskip

\proof[\bf{Proof of Theorem \ref{theorem2.11}}]
Let $f$ be any $\mu$-measurable function on $X$.
If $f$ does not belong to $L^{\alpha,\infty}(X)$ then $\left\|f\right\|^{\ast}_{\alpha,\infty}=\infty$ and there is nothing to prove. So we assume that $f$ is in $L^{\alpha,\infty}(X)$ and put $\left\|f\right\|^{\ast}_{\alpha,\infty}=A$.
\begin{enumerate}
\item[a)] Let us fix $r$ and $\lambda$ in $\left(0,\infty\right)$ and put
\begin{equation*}
E=\left\{x\in X:\left|f(x)\right|^{q}>\beta\right\}\text{ with }\beta=\frac{\lambda}{4\varphi(\rho^{m_{r}+1})\mathfrak b}.
\end{equation*}
For any integer $1\leq j<N_{m_{r}}$ such that $\left\|f\chi_{E^{m_{r}}_{j}}\right\|^{q}_{q}>\lambda$, we have
\begin{equation*}
\lambda-\left\|f\chi_{E\cap E^{m_{r}}_{j}}\right\|^{q}_{q}<\int_{E^{m_{r}}_{j}\setminus E}\left|f(x)\right|^{q}d\mu(x)\leq\beta\mu\left(E^{m_{r}}_{j}\setminus E\right)\leq\frac{\lambda}{4}.
\end{equation*}
Therefore $\frac{3\lambda}{4}<\left\|f\chi_{E\cap E^{m_{r}}_{j}}\right\|^{q}_{q}$ and
\begin{equation*}
\#\left(\left\{j:1\leq j<N_{m_{r}}\text{ and } \left\|f\chi_{E^{m_{r}}_{j}}\right\|^{q}_{q}>\lambda\right\}\right)\\
\leq\#\left(\left\{j:1\leq j<N_{m_{r}}\text{ and } \left\|f\chi_{E^{m_{r}}_{j}\cap E}\right\|^{q}_{q}>\frac{3\lambda}{4}\right\}\right).
\end{equation*}
Thus
\begin{equation*}
\begin{aligned}
&\frac{3\lambda}{4}\#\left(\left\{j:1\leq j<N_{m_{r}}\text{ and } \left\|f\chi_{E^{m_{r}}_{j}}\right\|^{q}_{q}>\lambda\right\}\right)\ \ \ \ \ \ \ \\ 
&\ \ \ \ \ \ \ \ \ \ \ \ \ \ \ \ \ \ \ \ \ \ \ \ \ \ \ \ \leq\frac{3\lambda}{4}\#\left(\left\{j:1\leq j<N_{m_{r}}\text{ and } \left\|f\chi_{E^{m_{r}}_{j}\cap E}\right\|^{q}_{q}>\frac{3\lambda}{4}\right\}\right)\\
&\ \ \ \ \ \ \ \ \ \ \ \ \ \ \ \ \ \ \ \ \ \ \ \ \ \ \ \ \leq\sum^{N_{m_{r}}}_{j=1}\int_{E\cap E^{m_{r}}_{j}}\left|f(x)\right|^{q}d\mu(x)\leq\left(\frac{\alpha}{\alpha-q}\right)A^{q}\mu\left(E\right)^{1-\frac{q}{\alpha}}
\end{aligned}
\end{equation*}
according to Kolmogorov condition (see \cite{GR}).
As 
\begin{equation*}
\mu\left(E\right)=\lambda_{f}(\beta^{\frac{1}{q}})\leq\left(\beta^{-\frac{1}{q}}A\right)^{\alpha}=\left(\frac{4\varphi(\rho^{m_{r}+1})\mathfrak b}{\lambda}\right)^{\frac{\alpha}{q}}A^{\alpha},
\end{equation*}
we obtain
\begin{equation*}
\frac{3\lambda}{4}\#\left(\left\{j:1\leq j<N_{m_{r}}\text{ and } \left\|f\chi_{E^{m_{r}}_{j}}\right\|^{q}_{q}>\lambda\right\}\right)\leq\frac{\alpha}{\alpha-q}\left(\frac{4\varphi(\rho^{m_{r}+1})\mathfrak b}{\lambda}\right)^{\frac{\alpha}{q}-1}A^{\alpha},
\end{equation*}
that is
\begin{equation}
\#\left(\left\{j:1\leq j<N_{m_{r}}\text{ and } \left\|f\chi_{E^{m_{r}}_{j}}\right\|^{q}_{q}>\lambda\right\}\right)\leq C \varphi(\rho^{m_{r}+1})^{\frac{\alpha}{q}-1}\lambda^{-\frac{\alpha}{q}}A^{\alpha}\label{85}
\end{equation}
with $C=\frac{4^{\frac{\alpha}{q}}\alpha\mathfrak b^{\frac{\alpha}{q}-1}}{3(\alpha-q)}.$
\item[b)] Assume that $p<\infty$. Suppose that $1<s<\infty$ and $r>0$ and put
\begin{equation*}
d_{j}=\left\|f\chi_{E^{m_{r}}_{j}}\right\|_{q}A^{-1}\left[\mathfrak b\varphi(\rho^{m_{r}+1})\right]^{\frac{1}{\alpha}-\frac{1}{q}}\left(\frac{\alpha}{\alpha-q}\right)^{-\frac{1}{q}},\ \ \ 1\leq j<N_{m_{r}}.
\end{equation*}
From Kolmogorov condition, we obtain $0\leq d_{j}\leq 1$ for $1\leq j<N_{m_{r}}$. In addition, for any number $\lambda$, we have 
\begin{equation*}
\#\left(\left\{j:1\leq j<N_{m_{r}},d_{j}>\lambda\right\}\right)\leq C\left[\left(\frac{\alpha}{\alpha-q}\right)^{\frac{1}{q}}\lambda\right]^{-\alpha}
\end{equation*}
according to Inequality (\ref{85}). 
Thus, we have
\begin{eqnarray*}
\sum^{N_{m_{r}}}_{j=1}d^{p}_{j}&=&\sum^{\infty}_{n=1}\left(\sum_{s^{-n-1}<d_{k}\leq s^{-n}}d^{p}_{k}\right)\leq\sum^{\infty}_{n=1}C\left[\left(\frac{\alpha}{\alpha-q}\right)^{\frac{1}{q}}s^{-n-1}\right]^{-\alpha}s^{-np}\\
&\leq&C\left(\frac{\alpha}{\alpha-q}\right)^{\frac{\alpha}{q}}\sum^{\infty}_{n=1}s^{\alpha-(p-\alpha)n}=C\left(\frac{\alpha}{\alpha-q}\right)^{\frac{\alpha}{q}}\frac{s^{2\alpha-p}}{s^{p-\alpha}-1}.
\end{eqnarray*}
This implies that
\begin{eqnarray*}
\left\|f\right\|^{dm_{r}}_{q,p,\alpha}&\leq&\left[\sup_{1\leq j<N_{m_{r}}}\frac{\varphi(\rho^{m_{r}+1})}{\mu\left(E^{m_{r}}_{j}\right)}\right]^{\frac{1}{q}-\frac{1}{\alpha}}\left\{\sum^{N_{m_{r}}}_{j=1}\left[\left\|f\chi_{E^{m_{r}}_{j}}\right\|_{q}A^{-1}\left(\mathfrak b\varphi(\rho^{m_{r}+1})\right)^{\frac{1}{\alpha}-\frac{1}{q}}\left(\frac{\alpha}{\alpha-q}\right)^{-\frac{1}{q}}\right]^{p}\right\}^{\frac{1}{p}}\\
&\times&A\mathfrak b^{\frac{1}{q}-\frac{1}{\alpha}}\left(\frac{\alpha}{\alpha-q}\right)^{\frac{1}{q}}\\
&\leq&\left[\sup_{1\leq j<N_{m_{r}}}\frac{\varphi(\rho^{m_{r}+1})}{\mu\left(E^{m_{r}}_{j}\right)}\right]^{\frac{1}{q}-\frac{1}{\alpha}}\left(\frac{\alpha}{\alpha-q}\right)^{\frac{1}{q}(1-\frac{\alpha}{p})}\left(\frac{s^{2\alpha-p}}{s^{p-\alpha}-1}\right)^{\frac{1}{p}}C^{\frac{1}{p}}A.
\end{eqnarray*}
As $r>0$ is arbitrary in $\left(0,\infty\right)$, we obtain
\begin{equation*}
\left\|f\right\|_{q,p,\alpha}\leq C\left\|f\right\|^{\ast}_{\alpha,\infty}
\end{equation*}
with $C$ a constant not depending on $f$.
\item[c)] For any number $r>0$ and positive integer $j<N_{m_{r}}$, we have according to Kolmogorov condition
\begin{equation*}
\mu\left(E^{m_{r}}_{j}\right)^{\frac{1}{\alpha}-\frac{1}{q}}\left\|f\chi_{E^{m_{r}}_{j}}\right\|_{q}\leq\left(\frac{\alpha}{\alpha-q}\right)^{\frac{1}{q}}A.
\end{equation*}
Thus
\begin{equation*}
\left\|f\right\|_{q,\infty,\alpha}\leq \left(\frac{\alpha}{\alpha-q}\right)^{\frac{1}{q}}\left\|f\right\|^{\ast}_{\alpha,\infty}.
\end{equation*}
\end{enumerate}
\epf

\medskip

Up to now we have used in our proofs the decomposition of X in dyadic cubes as given by Sawyer and Wheeden in \cite{sw}. The dyadic cubes $E^{k}_{j}$ $k\geq m,1\leq j<N_{k}$ have their size bounded below by $\rho^{m}$ with $\rho>1$ and $m$ a fixed integer. For the proof of the next theorem, we shall use the following decomposition given by Christ in \cite{C}.

\begin{lem}\label{christ} There exist a collection of open subsets $\left\{Q^{k}_{\alpha}\subset X:k\in\mathbb Z,\alpha\in I_{k}\right\}$, and  constants $\rho$ in $\left(0,1\right)$, $\mathfrak c_{0}>0$, $\eta>0$ and $\mathfrak c_{1},\mathfrak c_{2}<\infty$ such that
\begin{enumerate}
\item[(i)] $\mu\left(X\setminus\bigcup_{\alpha}Q^{k}_{\alpha}\right)=0\ \ \forall k,$
\item [(ii)]if $\ell\geq k$ then either $Q^{\ell}_{\beta}\subset Q^{k}_{\alpha}$ or $Q^{\ell}_{\beta}\cap Q^{k}_{\alpha}=\emptyset,$
\item [(iii)]for each $(k,\alpha)$ and each $\ell<k$ there is a unique $\beta$ such that $Q^{k}_{\alpha}\subset Q^{\ell}_{\beta},$
\item [(iv)]diameter$(Q^{k}_{\alpha})\leq \mathfrak c_{1}\rho^{k},$
\item [(v)]each $Q^{k}_{\alpha}$ contains some ball $B_{(z^{k}_{\alpha},\mathfrak c_{0}\rho^{k})}$,
\item [(vi)]$\mu\left(\left\{x\in Q^{k}_{\alpha}:d(x,X\setminus Q^{k}_{\alpha})\leq t\rho^{k}\right\}\right)\leq \mathfrak c_{2}t^{\eta}\mu(Q^{k}_{\alpha})\ \ \ \forall k,\alpha,\ \forall t>0$.
\end{enumerate} 
\end{lem}

\proof[\bf{Proof of Theorem \ref{theorem2.12}}]
Throughout the proof, we shall use the notation of the above lemma.
\begin{enumerate}
\item[A-]

\begin{enumerate}
\item Let us consider an element $\beta_{1}$ of $I_{1}$ and put $E_{1}=Q^{1}_{\beta_{1}}$. Then
\begin{equation*}
B_{(z^{1}_{\beta_{1}},\mathfrak c_{0}\rho)}\subset Q^{1}_{\beta_{1}}\subset B_{(z^{1}_{\beta_{1}},\mathfrak c_{1}\rho)}.
\end{equation*}
So by Inequalities (\ref{normale}), (\ref{Ip}) and (\ref{Ig}), 
\begin{equation*} 
\mu\left(E_{1}\right)=m\in\left[\mathfrak a\mathfrak a_{0}\mathfrak c^{D_{\mu}}_{0}\rho^{D_{\mu}},\mathfrak b\mathfrak b_{0}\mathfrak c^{\delta_{\mu}}_{1}\rho^{\delta_{\mu}}\right].
\end{equation*}
\item Let $\alpha_{2}\in I_{-2^{2}-1}$ such that $Q^{1}_{\beta_{1}}\subset Q^{-2^{2}-1}_{\alpha_{2}}.$

Put
\begin{equation*}
F_{1}=\emptyset,\ F_{2}=Q^{-2^{2}-1}_{\alpha_{2}}\text{ and } \tilde{J}_{2}=\left\{j\in I_{-2^{2}-1}:d(Q^{-2^{2}-1}_{j},F_{2})>\mathfrak c_{1}\rho^{-2^{2}-1}\right\}.
\end{equation*}
For each $j\in\tilde{J}_{2}$, let $\beta_{j}\in I_{2}$ be so that $d(z^{-2^{2}-1}_{j},Q^{2}_{\beta_{j}})<\rho^{2}.$
We have   
\begin{equation*}
\mu\left(Q^{2}_{\beta_{j}}\right)\in\left[\mathfrak a\mathfrak a_{0}\mathfrak c^{D_{\mu}}_{0}\rho^{2D_{\mu}},\mathfrak b\mathfrak b_{0}\mathfrak c^{\delta_{\mu}}_{1}\rho^{2\delta_{\mu}}\right],\ j\in\tilde{J}_{2}.
\end{equation*}
We can therefore choose a finite subset $J_{2}$ of $\tilde{J}_{2}$ such that
\begin{equation*}
\sum_{j\in J_{2}}\mu\left(Q^{2}_{\beta_{j}}\right)\in\left[m,m+\mathfrak b\mathfrak b_{0}\mathfrak c^{\delta_{\mu}}_{1}\right).
\end{equation*}
Let us take
\begin{equation*}
E_{2}=\cup_{j\in J_{2}}Q^{2}_{\beta_{j}}.
\end{equation*}
\item Let us consider for every $j\in J_{2}$ the element $\alpha_{j}$ of $I_{-2^{3}-1}$ such that $Q^{-2^{2}-1}_{j}\subset Q^{-2^{3}-1}_{\alpha_{j}}.$

Put
\begin{equation*}
F_{3}=\cup_{j\in J_{2}}Q^{-2^{3}-1}_{\alpha_{j}}\text{ and } \tilde{J}_{3}=\left\{j\in I_{-2^{3}-1}:d(Q^{-2^{3}-1}_{j},F_{3})>\mathfrak c_{1}\rho^{-2^{3}-1}\right\}.
\end{equation*}
For any $j\in\tilde{J}_{3}$, let $\beta_{j}\in I_{3}$ such that $d(z^{-2^{3}-1}_{j},Q^{3}_{\beta_{j}})<\rho^{3}$. We have 
\begin{equation*}
 \mu\left(Q^{3}_{\beta_{j}}\right)\in\left[\mathfrak a\mathfrak a_{0}\mathfrak c^{D_{\mu}}_{0}\rho^{3D_{\mu}},\mathfrak b\mathfrak b_{0}\mathfrak c^{\delta_{\mu}}_{1}\rho^{3\delta_{\mu}}\right],\ j\in \tilde{J}_{3}.
\end{equation*}
 Thus we can pick a finite subset $J_{3}$ in $\tilde{J}_{3}$ such that 
\begin{equation*}
\sum_{j\in J_{3}}\mu\left(Q^{3}_{\beta_{j}}\right)\in\left[m,m+\mathfrak b\mathfrak b_{0}\mathfrak c^{\delta_{\mu}}_{1}\rho^{3\delta_{\mu}}\right).
\end{equation*}
Put $E_{3}=\cup_{j\in J_{3}}Q^{3}_{\beta_{j}}.$
\item By iteration we obtain two sequences $(E_{n})_{n\geq 1}$ and $(F_{n})_{n\geq 1}$ such that

$\bullet\;\mu\left(E_{n}\right)\in\left[m,m+\mathfrak b\mathfrak b_{0}\mathfrak c^{\delta_{\mu}}_{1}\rho^{n\delta_{\mu}}\right)\text{ and } E_{n}=\cup_{j\in J_{n}}Q^{n}_{\beta_{j}}$, where $J_{n}$ is a finite subset of $I_{-2^{n}-1},$\\
$\bullet \ d(z^{-2^{n}-1}_{j},Q^{n}_{\beta_{j}})<\rho^{n}$ and $d(Q^{-2^{n}-1}_{j},F_{n})>\mathfrak c_{1}\rho^{-2^{n}-1}.$ 
\end{enumerate}
\item[B-] We fix $n\geq1$.
\begin{enumerate}
\item Let $(x,r)\in X\times\mathbb R^{\ast}_{+}.$ Suppose that $\ell,j\in J_{n}$ with $B_{(x,r)}\cap Q^{n}_{\beta_{j}}\neq\emptyset\neq B_{(x,r)}\cap Q^{n}_{\beta_{\ell}}$.

 There exists $x_{1},x_{2}\in Q^{n}_{\beta_{j}}$ and $y_{1},y_{2}\in Q^{n}_{\beta_{\ell}}$ such that $d(z^{-2^{n}-1}_{j},x_{1})<\rho^{n}$,  $x_{2}\in B_{(x,r)}$, $d(z^{-2^{n}-1}_{\ell},y_{1})<\rho^{n}$, $y_{2}\in B_{(x,r)}$.
 
  Therefore
\begin{eqnarray*}
\rho^{-2^{n}-1}&\leq&d\left(z^{-2^{n}-1}_{j},z^{-2^{n}-1}_{\ell}\right)\leq \kappa\left[d\left(z^{-2^{n}-1}_{j},x_{1}\right)+d\left(x_{1},z^{-2^{n}-1}_{\ell}\right)\right]\\
&\leq& \kappa\rho^{n}+2\kappa^{3}\mathfrak c_{1}\rho^{n}+2\kappa^{4}r+\kappa^{4}\left[d(y_{2},y_{1})+d(y_{1},z^{-2^{n}-1}_{\ell})\right]\\
&<&\left(\kappa+2\kappa^{3}\mathfrak c_{1}+2\kappa^{5}\mathfrak c_{1}+\kappa^{4}\right)\rho^{n}+2\kappa^{4}r.
\end{eqnarray*}
It follows that
\begin{eqnarray*}
r&>&\frac{1}{2\kappa^{4}}\left[\rho^{-2^{n}-1}-\left(\kappa+2\kappa^{3}\mathfrak c_{1}+2\kappa^{5}\mathfrak c_{1}\right)\rho^{n}\right]\\
&=&\frac{\rho^{n}}{2\kappa^{4}}\left[\rho^{-2^{n}-1-n}-\left(\kappa+2\kappa^{3}\mathfrak c_{1}+2\kappa^{5}\mathfrak c_{1}\right)\right].
\end{eqnarray*}
\item In the sequel we assume that $n$ is sufficiently great such that 
\begin{equation*}
 \mathfrak c_{1}\rho^{n}<1\leq\frac{\rho^{n}}{2\kappa^{4}}\left[\rho^{-2^{n}-1-n}-\left(\kappa+2\kappa^{3}\mathfrak c_{1}+2\kappa^{5}\mathfrak c_{1}+\kappa^{4}\right)\right]=r_{n}.
\end{equation*}

{\it $1^{rst}$ case.} We suppose that $0<r\leq \mathfrak c_{1}\rho^{n}$.

Then every ball $B_{(x,r)}$ meets at most one $Q^{n}_{\beta_{j}}$ $(j\in J_{n})$. Therefore,
\begin{eqnarray*}
_{r}\left\| \chi _{_{E_{n}}}\right\| _{q,p,\alpha}&=&\left[\int_{X}\left(\mu(B_{(x,r)})^{\frac{1}{\alpha}-\frac{1}{p}-\frac{1}{q}}\left\|\chi_{E_{n}\cap B_{(x,r)}}\right\|_{q}\right)^{p}d\mu(x)\right]^{\frac{1}{p}}\\
&=& \left[ \underset{j\in J_{n}}{\sum }\int_{\left\{x\in X:B_{(x,r)\cap Q^{n}_{\beta_{j}}\neq\emptyset}\right\}}\left(\mu(B_{(x,r)})^{\frac{1}{\alpha}-\frac{1}{p}-\frac{1}{q}}\mu\left(E_{n}\cap B_{(x,r)}\right)^{\frac{1}{q}}\right)^{p}d\mu(x)\right]^{\frac{1}{p}}\\
&\leq&\left[ \underset{j\in J_{n}}{\sum }\mu\left(\left\{x\in X:B_{(x,r)}\cap Q^{n}_{\beta_{j}}\neq\emptyset\right\}\right)\sup_{B_{(x,r)\cap Q^{n}_{\beta_{j}}\neq\emptyset}}\mu(B_{(x,r)})^{\frac{p}{\alpha}-1}\right]^{\frac{1}{p}}\\
&\leq&\left[ \underset{j\in J_{n}}{\sum }\mu\left(B_{(z^{n}_{\beta_{j}},\kappa(r+\mathfrak c_{1}\rho^{n}))}\right)(\mathfrak b\varphi(r))^{\frac{p}{\alpha}-1}\right]^{\frac{1}{p}}\\
&\leq&(\mathfrak b\varphi(r))^{\frac{1}{\alpha}-\frac{1}{p}}\left[ \underset{j\in J_{n}}{\sum }C_{\mu}(\frac{\kappa(r+\mathfrak c_{1}\rho^{n})}{\mathfrak c_{0}\rho^{n}})^{D_{\mu}}\mu(Q^{n}_{\beta_{j}})\right]^{\frac{1}{p}}\\
&\leq& C^{\frac{1}{p}}_{\mu}\left(\frac{2\kappa\mathfrak c_{1}}{\mathfrak c_{0}}\right)^{\frac{D_{\mu}}{p}}\mu\left(E_{n}\right)^{\frac{1}{p}}(\mathfrak b\mathfrak b_{0}\mathfrak c_{1}\rho^{n})^{(\frac{1}{\alpha}-\frac{1}{p})\delta_{\mu}}.
\end{eqnarray*}\label{parta}

{\it $2^{nd}$ case.} We suppose that $\mathfrak c_{1}\rho^{n}<r\leq r_{n}$.\\
 Arguing as in the first case,  we obtain
\begin{eqnarray*}
_{r}\left\| \chi _{_{E_{n}}}\right\| _{q,p,\alpha}&\leq&\left[ \underset{j\in J_{n}}{\sum }\mu\left(\left\{x\in X:B_{(x,r)}\cap Q^{n}_{\beta_{j}}\neq\emptyset\right\}\right)\sup_{B_{(x,r)}\cap Q^{n}_{\beta_{j}}\neq\emptyset}\mu(B_{(x,r)})^{\frac{p}{\alpha}-\frac{p}{q}-1}\mu\left(B_{(x,r)}\cap Q^{n}_{\beta_{j}}\right)^{\frac{p}{q}}\right]^{\frac{1}{p}}\\
&\leq& \left\{\underset{j\in J_{n}}{\sum }C_{\mu}(\frac{\kappa(r+\mathfrak c_{1}\rho^{n})}{r})^{D_{\mu}} \mathfrak b\varphi(r)\left(\mathfrak a\varphi(r)\right)^{\frac{p}{\alpha}-\frac{p}{q}-1}\mu\left( Q^{n}_{\beta_{j}}\right)\left[\mathfrak b\mathfrak b_{0}\left(\mathfrak c_{1}\rho^{n}\right)^{\delta_{\mu}}\right]^{\frac{p}{q}-1}\right\}^{\frac{1}{p}}.
\end{eqnarray*}
For the second inequality we have used the doubling condition of $\mu$, the relationship  between $\mu$ and $\varphi$, the growth condition on $\varphi$ and the inclusion $Q^{n}_{\beta_{j}}\subset B_{(z^{n}_{\beta_{j}},\mathfrak c_{1}\rho^{n})}$.
Thus
\begin{eqnarray*}
{r}\left\| \chi _{_{E_{n}}}\right\| _{q,p,\alpha}&\leq& C\varphi(r)^{\frac{1}{\alpha}-\frac{1}{q}}\mu\left( E_{n}\right)^{\frac{1}{p}}\left(\rho^{n}\right)^{\delta_{\mu}(\frac{1}{q}-\frac{1}{p})}\\
&\leq& C\mu\left(E_{n}\right)^{\frac{1}{p}}\varphi\left(\rho^{n}\right)^{\frac{1}{\alpha}-\frac{1}{q}}\left(\mathfrak c_{1}\rho^{n}\right)^{\delta_{\mu}(\frac{1}{q}-\frac{1}{p})}
\leq C\mu\left( E_{n}\right)^{\frac{1}{p}}\left(\mathfrak c_{1}\rho^{n}\right)^{\delta_{\mu}(\frac{1}{q}-\frac{1}{p})+D_{\mu}(\frac{1}{\alpha}-\frac{1}{q})}.
\end{eqnarray*}

{\it $3^{rd}$ case.} We suppose $r>r_{n}$.
\begin{eqnarray*}
_{r}\left\| \chi _{_{E_{n}}}\right\| _{q,p,\alpha}
&\leq& \left[ \underset{j\in J_{n}}{\sum }\int_{\left\{x\in X:B_{(x,r)}\cap Q^{n}_{\beta_{j}}\neq\emptyset\right\}}\mu(B_{(x,r)})^{\frac{p}{\alpha}-\frac{p}{q}-1}\mu\left(E_{n}\cap B_{(x,r)}\right)^{\frac{p}{q}}d\mu(x)\right]^{\frac{1}{p}}\\
&\leq& \left[ \mu\left(E_{n}\right)^{\frac{p}{q}}\underset{j\in J_{n}}{\sum }\mu(B_{(z^{n}_{\beta_{j}},\kappa(r+\mathfrak c_{1}\rho^{n}))})(\mathfrak a\varphi(r))^{\frac{p}{\alpha}-\frac{p}{q}-1}\right]^{\frac{1}{p}}\\
&\leq&  \mu\left(E_{n}\right)^{\frac{1}{p}}(\mathfrak a\varphi(r))^{\frac{1}{\alpha}-\frac{1}{q}-\frac{1}{p}}\left[\#(J_{n})C_{\mu}(2\kappa)^{D_{\mu}}\varphi(r)\right]^{\frac{1}{p}}.
\end{eqnarray*}
But for all $j\in J_{n}$, $\mathfrak a_{0}\mathfrak c^{D_{\mu}}_{0}\rho^{nD_{\mu}}\leq\mu(Q^{n}_{\beta_{j}})\leq b_{0}\mathfrak c_{1}^{\delta_{\mu}}\rho^{n\delta_{\mu}}$. Thus 
\begin{equation*}
\#(J_{n})\leq\frac{m}{\mathfrak a_{0}\mathfrak c^{D_{\mu}}_{0}\rho^{nD_{\mu}}},
\end{equation*}
and 
\begin{eqnarray*}
_{r}\left\| \chi _{_{E_{n}}}\right\| _{q,p,\alpha}&\leq&C\mu\left(E_{n}\right)^{\frac{1}{p}}\varphi(r)^{\frac{1}{\alpha}-\frac{1}{q}}\frac{1)^{\frac{1}{p}}}{\rho^{nD_{\mu}/p}}\leq C\mu\left(E_{n}\right)^{\frac{1}{p}}r^{\delta(\frac{1}{\alpha}-\frac{1}{q})}_{n}\rho^{-nD_{\mu}/p}\\
&\leq& C\mu\left(E_{n}\right)^{\frac{1}{p}}\left(\rho^{-n}\right)^{\frac{D_{\mu}}{p}+\delta_{\mu}(\frac{1}{q}-\frac{1}{\alpha})}\left[\rho^{-2^{n}-1-n}-(\kappa+2\kappa^{3}\mathfrak c_{1}+2\kappa^{5}\mathfrak c_{1})\right]^{\delta_{\mu}(\frac{1}{\alpha}-\frac{1}{q})}\\
&\leq& C\mu\left(E_{n}\right)^{\frac{1}{p}}\left(\rho^{-n}\right)^{D_{\mu}(\frac{1}{p}-\frac{1}{q})+\delta(\frac{1}{q}-\frac{1}{p})}\frac{(\rho^{-n})^{D_{\mu}(\frac{1}{\alpha}-\frac{1}{q}+\frac{1}{p})+\delta_{\mu}(\frac{2}{q}-\frac{1}{p}-\frac{1}{\alpha})}}{(\rho^{-2^{n}-1-n}-(\kappa+2\kappa^{3}\mathfrak c_{1}+2\kappa^{5}\mathfrak c_{1}))^{\delta_{\mu}(\frac{1}{q}-\frac{1}{\alpha})}}.
\end{eqnarray*}
It follows that if we choose $n_{0}$ such that for all $n\geq n_{0}$
\begin{equation*}
\frac{(\rho^{-n})^{D_{\mu}(\frac{1}{\alpha}-\frac{1}{q}+\frac{1}{p})+\delta_{\mu}(\frac{2}{q}-\frac{1}{p}-\frac{1}{\alpha})}}{(\rho^{-2^{n}-1-n}-(\kappa+2\kappa^{3}\mathfrak c_{1}+2\kappa^{5}\mathfrak c_{1}))^{\delta_{\mu}(\frac{1}{q}-\frac{1}{\alpha})}}<1,
\end{equation*}
then  
\begin{equation*}
_{r}\left\| \chi _{_{E_{n}}}\right\| _{q,p,\alpha}\leq C\rho^{-n})^{D_{\mu}(\frac{1}{\alpha}-\frac{1}{q})+\delta_{\mu}(\frac{1}{q}-\frac{1}{p})}\mu\left(E_{n}\right)^{\frac{1}{p}},\ \ \ r>0.
\end{equation*}
That is 
\begin{equation*}
\left\| \chi _{_{E_{n}}}\right\| _{q,p,\alpha}\leq C\rho^{-n})^{D_{\mu}(\frac{1}{\alpha}-\frac{1}{q})+\delta_{\mu}(\frac{1}{q}-\frac{1}{p})}
\end{equation*}
and therefore
\begin{equation*}
\left\|\chi_{\cup_{n\geq n_{0}}E_{n}}\right\|_{q,p,\alpha}<\infty\text{ and }\left\|\chi_{\cup_{n\geq n_{0}}E_{n}}\right\|^{\ast}_{\alpha,\infty}=\infty.
\end{equation*}
 \end{enumerate}
\end{enumerate}
Thus, $f=\chi_{\cup_{n>n_{0}}E_{n}}$ is in $(L^{q},L^{p})^{\alpha}(X)\setminus L^{\alpha,\infty}(X)$.
\epf

\end{document}